\documentclass[onefignum,onetabnum]{siamart190516}

\usepackage{lipsum}
\usepackage{amsfonts}
\usepackage{graphicx}
\usepackage{epstopdf}
\usepackage{algorithmic}
\usepackage{bbm}  
\usepackage{multirow}
\ifpdf
  \DeclareGraphicsExtensions{.eps,.pdf,.png,.jpg}
\else
  \DeclareGraphicsExtensions{.eps}
\fi

\newtheorem{sketch}{Sketch}

\newcommand{\RR}{\mbox{{\sl I}}\!\mbox{{\sl R}}} 

\newcommand{\ds}{\displaystyle}

\newcommand{\Inn}{\textsf{Inner}}
\newcommand{\Out}{\textsf{Outer}}

\newsiamremark{remark}{Remark}
\newsiamremark{hypothesis}{Hypothesis}
\crefname{hypothesis}{Hypothesis}{Hypotheses}
\newsiamthm{claim}{Claim}

\headers{Solving high dimensional SABR/LIBOR PDE models}{J.G. L\'opez-Salas, S. P\'erez-Rodr\'iguez, and C. V\'azquez}

\title{AMFR-W numerical methods for solving high dimensional SABR/LIBOR PDE models\thanks{Submitted to the editors DATE.
\funding{ This work was funded by the Spanish Grants MTM2016-76497-R, PID2019-108584-
RB-I00 and MTM2016-77735-C3-3-P, as well as Xunta de Galicia grant ED431C2018/033, all includ-
ing FEDER funding. First and third authors acknowledge the support received from the Centro de
Investigacin de Galicia ”CITIC”, funded by Xunta de Galicia and the European Union (European
Regional Development Fund- Galicia 2014-2020 Program), by grant ED431G 2019/01.}}}

\author{J.G. L\'opez-Salas\thanks{Department of Mathematics and CITIC, Universidade da Coru\~na, Campus de Elvi\~na s/n, 15071 - A Coru\~na, Spain (\email{jose.lsalas@udc.es}).}
\and S. P\'erez-Rodr\'iguez\thanks{Departamento de An\'alisis Matem\'atico, Universidad de La Laguna, 38208, La Laguna, Tenerife,
Canary Islands, Spain
  (\email{sperezr@ull.edu.es}).}
\and C. V\'azquez\thanks{Department of Mathematics and CITIC, Universidade da Coru\~na, Campus de Elvi\~na s/n, 15071 - A Coru\~na, Spain (\email{carlosv@udc.es}).}}

\usepackage{amsopn}

\ifpdf
\hypersetup{
  pdftitle={Solving high dimensional SABR/LIBOR PDE models},
  pdfauthor={J.G. L\'opez-Salas, S. P\'erez-Rodr\'iguez, and C. V\'azquez}
}
\fi

\begin{document}

\maketitle

\begin{abstract}
In this work we mainly develop a new numerical methodology to solve a PDE model
recently proposed in the literature for pricing interest rate derivatives. More precisely, we use high
order in time AMFR-W methods, which belong to a class of W-methods based on Approximate
Matrix Factorization (AMF) and are specially suitable in the presence of mixed spatial derivatives.
High order convergence in time allows larger time steps which combined with the splitting of the
involved operators, highly reduces the computational time for a given accuracy. Moreover, the con-
sideration of a large number of underlying forward rates makes the PDE problem high dimensional in
space, so the use of AMFR-W methods with a sparse grids combination technique represents another
innovative aspect, making AMFR-W more efficient than with full grids and opening the possibility of
parallelization. Also the consideration of new homogeneous Neumann boundary conditions provides
another original feature to avoid the difficulties associated to the presence of boundary layers when
using Dirichlet ones, specially in advection dominated regimes. These Neumann boundary condi-
tions motivate the introduction of a modified combination technique to overcome a decrease in the
accuracy of the standard combination technique.
\end{abstract}

\begin{keywords}
SABR-LIBOR market models, high dimensional PDEs, AMFR-W methods, finite differences,  sparse grids combination technique
\end{keywords}

\begin{AMS}
  65M06, 65M20, 65M50, 65W10, 91G30, 91G80
\end{AMS}

\section{Introduction}

High dimensional parabolic Partial Differential Equations \\ (PDEs) arise in many fields of science and engineering problems, as for example in computational biology for stochastic gene networks \cite{abpv18} or in computational finance for pricing financial derivatives \cite{lopezvazquez18}, when a large number of underlying stochastic factors are involved in their equivalent stochastic formulations. In fact, each stochastic factor gives rise to one spatial-like variable in the corresponding PDE. In this high dimensional setting, when using finite differences for the spatial discretization, the complexity of standard grid based approaches grows exponentially with the dimensions of the problem as well as the computational times, thus giving rise to the so called {\em curse of dimensionality}. Thus, alternative techniques to the standard full grid are required. Also, the use of high order time integration schemes turns out to be very convenient to allow larger time steps and therefore reduce the computational time to get a prescribed accuracy.

In the present work we mainly propose a new numerical technique for solving the high dimensional PDE problem governing the Stochastic Alpha Beta Rho-LIBOR
Market Model (SABR-LMM) PDE model introduced in \cite{lopezvazquez18}. The more classical LMM has been introduced for pricing interest rate derivatives which depend on the evolution of a certain number of forward LIBOR rates, when their volatility was assumed to be constant (see \cite{brigo_mercurio2007}, for example). More recently, the consideration of stochastic volatility has been incorporated for a better fit to market data by combining the classical LMM and the SABR model for stochastic volatility in \cite{haganlesniewsky08, mercuriomorini09, rebonatowhite09}, with different modelling approaches. In these works, the number of stochastic factors depends on the number of forward rates and volatilities that are considered. Their formulations are posed in terms of expectations to be solved by means of Monte Carlo techniques. Although we are aware of the evolution of LMM as a consequence of the financial crisis in 2007 to incorporate the practical presence of a multicurve setting (see \cite{mercurio12}, for example) and the recent ongoing studies related to LIBOR transitions not consolidated in market practice yet \cite{lyashenkomercurio20}, we have chosen to start from the classical version of LMM as described in \cite{brigo_mercurio2007}.

More recently, in \cite{lopezvazquez18} a PDE formulation is obtained for the Mercurio and Morini model in \cite{mercuriomorini09} and a set of numerical methods are proposed to solve it. More precisely, the combination of standard finite differences in space and a $\theta-$method in time are proposed on uniform full grids. Moreover, by arguing that these standard finite difference methods based on traditional full grid are not able to price interest rate derivatives with more than three or four stochastic forward rates, a sparse grid combination technique is applied. A rigorous analysis of finite differences schemes
in the sparse grid combination technique in arbitrary dimensions is carried out in \cite{ReisingerWittum2007}.In order to overcome the curse of dimensionality one can try to use high order discretizations in time and space, mainly to reduce the required discretization points to achieve certain accuracy. Although for sparse grids in space there are lots of works (see \cite{hendricksHeuerEhrhardtGunther2017} and the references therein), for time discretization only schemes up to order two have been applied. Following this objective, unlike in \cite{lopezvazquez18}, in the present article we propose higher order discretization techniques in time based on a special class of W-methods \cite{steihaug79aat}, the AMFR-W methods introduced in \cite{amfrw18}.

These AMFR-W-methods are specially suitable for parabolic problems involving mixed spatial derivatives as it is the case in the SABR-LMM model proposed in \cite{lopezvazquez18}. In \cite{amfrw18} it is shown that they are unconditionally stable regardless the spatial dimension on linear constant coefficients PDEs with mixed derivative terms with both periodic boundary conditions and homogeneous Dirichlet boundary conditions. Moreover,  as the ADI methods \cite{KarelMishra2013,hendricksHeuerEhrhardtGunther2017}, the AMFR-W scheme takes advantage of the structure of the linear system obtained from the spatial discretization, so it can be decomposed into tridiagonal systems which can be solved in linear run-time. Thus, the computational effort is significantly reduced. However, while classical ADI schemes are of order two in time, the proposed AMFR-W scheme exhibits order three in time when standard full grids are used. An additional innovative aspect of the present work is the application of the AMFR-W methods in the context of sparse grids combination technique, which turns out to be an efficient tool for solving the SABR-LMM model in the required high dimensional setting. Furthermore, the introduction of more appropriate Neumann boundary conditions motivates the consideration of a modified combination technique to improve the convergence.

LMMs are usually simulated by means of Monte Carlo method, in contrast with the here proposed PDE methodology. We aim to avoid the slow Monte Carlo rate of convergence, $O(1/\sqrt{M})$ for all dimensions, $M$ being the number of simulations. 

The plan of the article is the following. In Section 2 we present the PDE model and justify the introduction of new homogeneous Neumann boundary conditions at the so called outflow boundaries. In Section 3, we introduce the space discretization of the PDE problem with finite differences to obtain an ODE system in suitable form for the application of the AMFR-W method. In Section 4 we describe the application of AMFR-W method to the ODE system to get the fully discretized problem. Section 5 is devoted to the methodology of sparse grids, including the standard and the modified combination techniques. In Section 6 we present and discuss the numerical results obtained for full and sparse grids. Finally, Section 7  contains some conclusions.


\section{PDE formulation of the SABR-LMM model}

As indicated in the previous section, we mainly address a new and more efficient numerical solution of the PDE formulation introduced in  \cite{lopezvazquez18} for the SABR-LMM proposed by Mercurio and Morini in \cite{mercuriomorini09} to price a financial derivative which depends on a certain number of forward LIBOR rates, i.e a swaption. In this section we introduce the PDE model and we incorporate some new boundary conditions. Although we address the reader to \cite{lopezvazquez18} for the statement of the model and further details, we need to introduce some financial concepts and their notations related to interest rates derivatives to be used along this article. In this respect, we also address the reader to the textbook \cite{brigo_mercurio2007}.

A zero coupon \textit{bond} with maturity at time $T$ pays its holder one unit of currency at time $T$. The zero coupon value at time $t<T$  will be denoted by $P(t,T)$, and is also referred as the discount factor from time $T$ to time $t$. A \textit{tenor structure} is defined as a set of ordered payment dates $T_0<T_1<\ldots<T_{N-1}<T_N$. The time gap between two consecutive tenor dates is denoted by $\tau_i = T_{i+1}-T_i$. In view of previous definitions, a payment of $x$ units at time $T_i$ is worth $xP(t,T_i)$ at time $t<T_i$.

 Next, we introduce the forward rates that enter in the LMM. We consider the \textit{forward} interest rate $F_i(t)$ as an interest rate we can contract to borrow or lend money during the future time period $[T_i,T_{i+1}]$, which is fixed at time $T_i$. Moreover, the value of $F_i(t)$ can be expressed in terms of discount factors in the form: $$F_i(t) = F(t;T_i,T_{i+1}) = \dfrac{1}{\tau_i}\left(\dfrac{P(t,T_i)}{P(t,T_{i+1})}-1\right) \mbox{ where } t\leq T_i.$$

Conversely, the price of a zero coupon bond at time $T_i$ that matures at $T_j$, $P(T_i,T_j)$, can be expressed in terms of forward LIBOR rates as $$P(T_i,T_j) = \displaystyle\prod_{k=i}^{j-1}\dfrac{1}{1+\tau_k F_k(T_i)}.$$

Among all interest rate derivatives, the simplest one is the caplet. A \textit{caplet} is a European call option on a forward rate. Thus, if the maturity of a caplet is $T_{i+1}$, at that time the holder of the caplet receives the payoff $\tau_i (F_i(T_i)-K)^+$, so its discounted payoff at time $t<T_{i+1}$ is given by
$P(t,T_{i+1})\tau_i (F_i(T_i)-K)^+,$ where $(\cdot)^+$ denotes the function $\max(\cdot,0)$ and $K$ is the strike (a fixed interest rate) of the caplet. If constant volatilities are assumed as it is the case in the classical LMM, the caplet price can be analytically computed with a Black's formula (see \cite{brigo_mercurio2007}, for details).

An interest rate \textit{swap} (IRS) is a contract to exchange interest payments at future fixed dates. At every time instant in the prescribed set of dates $T_{a+1},\ldots,T_{b}$ the contract holder pays a fixed interest rate $K$ and receives a floating forward LIBOR rate $F_i(T_i)$, which is fixed at time $T_i$. At time $T_a$ the value of the  IRS is given by
\begin{equation}
\mbox{IRS}(T_a;T_a,\ldots,T_b) = \displaystyle\sum_{i=a}^{b-1} P(T_a,T_{i+1})\tau_i (F_i(T_a)-K). \label{eq:swap}
\end{equation}

A European $T_a \times (T_b-T_a)$ \textit{swaption} is an option that gives the right to enter a swap at the future time $T_a$ (swaption maturity). The underlying swap length $T_b-T_a$ is referred as the tenor of the swaption. Therefore, the discounted swaption payoff to time $t$ is equal to $P(t,T_a)\big(\mbox{IRS}(T_a;T_a,\ldots,T_b)\big)^+.$

In the forthcoming section devoted to numerical results, several examples address the pricing of caplets and swaptions. Note that the payoff of a caplet just involves one forward rate, therefore its price at any time before maturity only depends on this particular forward rate. In the case of IRS or swaptions, their payoffs at expiration date depend on a certain number of forward rates, so their price at any time before expiration will also depend on them.

In Mercurio and Morini model it is assumed that a generic European interest rate derivative depends on the evolution of $N-1$ forward rates, $F_1,\, F_2,\dots,\, F_{N-1}$ associated to bonds with maturity related to the tenor structure, and a common stochastic volatility, $V$. Let $u=u(t,F_1,F_2,\dots,F_{N-1},V)$ denote the value at time $t\in [0,T]$ of this European interest rate derivative, with forward rates $F_i\in [0,F_i^{max}],\,i=1,2,\dots,N-1,$ and volatility $V\in [0,V^{max}]$. Note that the previously described caplets, IRS and swaptions are particular cases. As stated in \cite{lopezvazquez18}, the function $u$ satisfies the following PDE
\begin{equation}\label{mmNd_no}
\begin{array}{rl}
\ds\frac{\partial u}{\partial t} +
\ds\frac{1}{2}\sigma^2 V^2 \ds\frac{\partial^2 u}{\partial V^2} +
\ds\frac{1}{2}V^2  \ds\sum_{i,j=1}^{N-1}   \alpha_i\alpha_j \rho_{ij} F_i^{\beta} F_j^{\beta} \ds\frac{\partial^2 u}{\partial F_i\partial F_j}
& \\[0.7pc]
+\sigma V^2 \ds\sum_{i=1}^{N-1}   \alpha_i\phi_i F_i^{\beta} \ds\frac{\partial^2 u}{\partial F_i\partial V}
+
\ds\sum_{i=1}^{N-1} \mu_i F_i^{\beta} \ds\frac{\partial u}{\partial F_i}
& =0,
\end{array}
\end{equation}
where $\sigma$ is the volatility of the stochastic volatility $V$, $\alpha_i$ is a deterministic (constant) instantaneous volatility coefficient of the forward rate $F_i$, $\rho_{ij}$ is the correlation between the forward rates $F_i$ and $F_j$, $\phi_i$ is the correlation between $F_i$ and the stochastic volatility $V$ and $\mu_i$ is the drift of the $i$-th forward rate. Moreover, when the bond $P(t,T_1)$ is chosen as numeraire, the drifts $\mu_i$ depend on the forward rates as 
$\mu_1=0$, $\mu_i=\alpha_i V^2  \sum_{j=2}^i \frac{\tau_j F_j^{\beta}}{1+\tau_j F_j} \rho_{ij} \alpha_j$, $i\ge 2$.
 The parameter $\beta \in [0,1]$ is the elasticity of variance, which usually is 0, 0.5 or 1,  that corresponds to stochastic volatilities with normal, CIR or log-normal dynamics.
 For the correlation structure, as in \cite{lopezvazquez18} we consider the expression $\rho_{ij} = \mbox{e}^{-\lambda|T_i-T_j|}$, 
which depends on the constant parameter $\lambda$.

In view of the form of the differential operator governing the PDE, (\ref{mmNd_no}) must be completed with a final condition $u(T,F_1,F_2,\dots,F_{N-1},V)= g(T,F_1,F_2,\dots,F_{N-1})$, where the $g$ represents the derivative payoff, the expression of which depends on the interest rate derivative we are dealing with.

%

 In order to apply the method of lines (MoL) to discretize the previous model it is more convenient to write \eqref{mmNd_no} in terms of the time to maturity $T-t$ instead of the physical time $t$, so that the final condition turns into an initial condition. In an abuse of notation we keep the notation $t$ for the new formulation after this change in the time variable. More precisely, we rewrite the model by denoting $u=u(t,F_1,F_2,\dots,F_{N-1},F_N)$ the value of the interest rate derivative at time $T-t$, $t\in [0,T]$,  with forward rates $F_i\in [0,F_i^{max}],\,i=1,2,\dots,N-1,$ and volatility $V\in [0,V^{max}]$. After some easy calculus, we obtain the equation
\begin{equation}\label{mmNd}
\ds\frac{\partial u}{\partial t} =
 \ds\sum_{i=1}^{N}
d_i \ds\frac{\partial^2 u}{\partial F_i^2} +
\ds\sum_{i=1}^{N-1} \ds\sum_{k=i+1}^{N} m_{ik}\ds\frac{\partial^2 u}{\partial F_i\partial F_k}+
\ds\sum_{i=2}^{N-1} a_i \ds\frac{\partial u}{\partial F_i},
\end{equation}
with the {\it initial condition}
\begin{equation}\label{icNd}
u(0,F_1,F_2,\dots,F_{N-1},F_N)= g(T,F_1,F_2,\dots,F_{N-1}),
\end{equation}
where:
\begin{equation}\label{coef_diff}
d_i=d_i(F_i,F_N)= \left\{ \begin{array}{ll}
\ds\frac{1}{2}\alpha_i^2 \rho_{ii} F_i^{2\beta} F_N^2, & \hbox{if } 1\le i\le N-1, \\[0.5pc]
\ds\frac{1}{2}\sigma^2 F_N^2, & \hbox{if } i=N. \end{array}\right.
\end{equation}
Moreover, for $i=1,\dots,N-1$, $k=i+1,\dots,N$, we define
\begin{equation}\label{coef_mix}
m_{ik}=m_{ik}(F_i,F_k,F_N)= \left\{ \begin{array}{ll}
\alpha_i \alpha_k \rho_{ik} F_i^{\beta} F_k^{\beta} F_N^2, & \hbox{if } i+1\le k\le N-1, \\[0.5pc]
\alpha_i \sigma \phi_i F_i^{\beta} F_N^2, & \hbox{if } k=N. \end{array}\right.
\end{equation}
Note that the last term in \eqref{mmNd} is only defined for $N\ge 3$. In this case, for each $i=2,\dots,N-1$ we define:
\begin{equation}\label{coef_adv}
a_{i}=a_{i}(F_2,\dots,F_i,F_N)= \left(
\ds\sum_{j=2}^i  \alpha_i \alpha_j \rho_{ij} \Phi_{\beta}( F_j, \tau_j) \right) F_i^{\beta} F_N^2,
\end{equation}
where  $\Phi_{\beta}$ is the scalar function $\Phi_{\beta}(x,\tau):=\ds\frac{\tau x^{\beta}}{1+\tau x}$, $\tau>0$, $x\ge 0$.

In next paragraphs we will discuss about the appropriate boundary conditions to add to (\ref{mmNd})-(\ref{icNd}) to define the initial-boundary value problem. First, note that in \cite{lopezvazquez18} the following time-independent boundary conditions were considered to complete the formulation (\ref{mmNd})-(\ref{icNd}):
\begin{equation}\label{bound1_0}
\begin{array}{l}
\hbox{if } F_j=0 \, \hbox{ or}\, F_j=F_j^{max},\,1\le j\le N-1,   \\[0.6pc]
\hspace{2cm}u(t,F_1,\dots,F_{N-1},F_N)=g(T,F_1,\dots,F_{N-1}),\\[0.6pc]
u(t,F_1,\dots,F_{N-1},0)=u(0,F_1,\dots,F_{N-1},0)=g(T,F_1,\dots,F_{N-1}),   \\[0.5pc]  \ds\frac{\partial u}{\partial V}(t,F_1,\dots,F_{N-1},V^{max})=0.
\end{array}
\end{equation}

These boundary conditions are appropriate when $\beta=0$ since the PDE coefficients are independent of the forward rates ($F_j,\,j\le N-1$) and the advection terms are moderate. When $\beta>0$, they also can be appropriate when we have a small number of forward rates and a derivative with short maturity, but as soon as these numbers increase and therefore the dimension of the PDE, the advection-dominance of the PDE in \eqref{mmNd} becomes more relevant. It is well-known that imposing a Dirichlet boundary condition at the outflow boundaries $F_i=F_i^{max},\,i=1,\dots,N_1,$ in an advection-dominated setting could originate boundary layers (see, for instance, \cite[Sect. I.5]{hundsdorfer03nso}), and some previous numerical results with this model corroborate that.

In order to show that when $\beta \in (0,1]$ the larger the dimension $N$  the more advection-dominant the PDE \eqref{mmNd} becomes, it is better to write the PDE \eqref{mmNd} in the following conservative form
\begin{equation}\label{adr}
u_t+ \nabla \cdot (\underline{c}u)= \nabla \cdot (D\nabla u) + s(\underline{x},t,u),\quad \underline{x} \in \Omega \subset \mathbb{R}^N,\,t>0 \, ,
\end{equation}
where $\underline{x}=(F_1,\dots,F_N)$, $\underline{c}=(c_1,\dots,c_N)^T$, $D=\left(D_{ik}\right)_{i,k=1}^N$, $c_i=c_i(\underline{x},t)$, $D_{ik}=D_{ik}(\underline{x},t)$.

%
Note that in multi-dimensional advection-diffusion-reaction PDEs with variable coefficients, some additional advection terms arise in the conservative form coming from the partial derivatives of the coefficients of the second order diffusion terms.
%
%

After some manipulations it can be proved that
PDE \eqref{mmNd} admits an expression of type \eqref{adr} with $D_{ii}=d_i,\,D_{ik}=m_{ik}/2,\,\forall i\ne k  $, and
\begin{equation}\label{cis}
\begin{array}{l}c_i= \left(P_i^+ -P_i^- \right) \alpha_i F_i^{\beta} F_N^2 + \alpha_i \phi_i  F_i^{\beta} \sigma F_N, \quad 1\le i\le N-1, \\[0.5pc] c_N= P_N^+  \sigma F_N^2 +  \sigma^2 F_N,\end{array}
\end{equation}
where
$$\begin{array}{rl}
P_i^+= &\alpha_i \ds\frac{\beta}{F_i^{1-\beta}} + \ds\sum_{{\scriptsize \begin{array}{c} j=1 \\ j\ne i \end{array}}}^{N-1} \ds\frac{1}{2} \alpha_j \rho_{ij} \ds\frac{\beta}{F_j^{1-\beta}} , \quad 1\le i\le N-1, \\[0.5pc]
P_N^+= & \ds\sum_{ j=1}^{N-1} \ds\frac{1}{2} \alpha_j \phi_{j} \ds\frac{\beta}{F_j^{1-\beta}}, \quad P_1^-= 0, \quad P_2^-=\alpha_2 \Phi_{\beta}(F_2,\tau_2) \\[0.7pc]
 P_i^-=&\alpha_i \Phi_{\beta}(F_i,\tau_i) +\ds\sum_{j=2}^{i-1}\alpha_j \rho_{ij} \Phi_{\beta}(F_j,\tau_j)
, \quad 3\le i\le N-1.
\end{array}$$
%
%
Moreover, the reaction term takes the form $s(\underline{x},t,u)= \delta u$, $\delta =P\, F_N^2+Q\,  \sigma F_N+\sigma^2$, where $P$ and $Q$ depend on the forward rates $F_1,\dots,F_{N-1}$.

Obviously, when $0<\beta<1$ the PDE \eqref{mmNd} is equivalent to \eqref{adr} whenever $F_j>0,\,\forall j=1,\dots,N-1$, due to the lack of differentiability of $\underline{c}$ when $F_j=0$. Since we intend to use this expression \eqref{adr} to study the behaviour of the PDE at the boundaries $F_j=F_j^{max}$, we will assume for now that $F_j>0,\,  1\le j\le N-1$. On the other hand, when $\beta=1$ this additional assumption is not needed.

%
In most of practical cases we have that
\begin{equation}\label{requer1}
\tau_i F_i^{max} \le \ds\frac{\beta}{2-\beta},\quad 2\le i\le  N-2,\qquad \tau_{N-1} F_{N-1}^{max} \le \ds\frac{\beta}{1-\beta},
\end{equation}
what implies that, for $2\le i\le N-1$\footnote{In the sequel,  $\sum_{i=j}^{k} (\cdot)=0$ when $j>k$.},
\begin{equation}\label{pes}
\begin{array}{rl}
P_i^+ - P_i^-=&\alpha_i \left(\ds\frac{\beta}{F_i^{1-\beta}}- \Phi_{\beta}(F_i,\tau_i) \right)+ \ds\sum_{j=2}^{i-1} \frac{1}{2} \alpha_j \rho_{ij} \left(\ds\frac{\beta}{F_j^{1-\beta}}-2 \Phi_{\beta}(F_i,\tau_i) \right)  \\[0.5pc] &+ \ds\frac{1}{2} \alpha_1 \rho_{i1} \ds\frac{\beta}{F_1^{1-\beta}}+   \ds\sum_{j=i+1}^{N-1} \frac{1}{2} \alpha_j \rho_{ij}\ds\frac{\beta}{F_j^{1-\beta}} \ge 0.
\end{array}
\end{equation}
We must observe that in the case $\beta=1$,
$$\begin{array}{l}
\ds\frac{\beta}{F_i^{1-\beta}}- \Phi_{\beta}(F_i,\tau_i) =1-\ds\frac{\tau_iF_i}{1+\tau_i F_i}=\frac{1}{1+\tau_i F_i} \ge 0, \\[0.7pc]
\ds\frac{\beta}{F_j^{1-\beta}}- 2\Phi_{\beta}(F_j,\tau_j) =1-2\ds\frac{\tau_iF_i}{1+\tau_i F_i}=\frac{1-\tau F_j}{1+\tau_i F_j} \ge 0\quad  \hbox{ if } \tau_j F_j \le 1,
\end{array}$$
and the second requirement in \eqref{requer1} is superfluous.

Therefore, $c_i\ge 0$, $i=1,\dots,N$, so it is clear that all the boundaries $F_i=F_i^{max}$ $ 1\leq i\leq N$ are outflow boundaries since the outward normal vector on each one of these borders is the $N-$dimensional canonical vector $\mathbf{e}_i=(e_{ik})_{k=1}^N$, $e_{ii}=1, \, e_{ik}\ne 0, \,\forall k\ne i$.

On the other hand, we need to take into account that these advection coefficients depend on $N$,  $c_i=c_i^{(N)}$, and  $F_N=V$ so if we increase the number of forward rates from $\{F_1,\dots,F_{N-2}\}$ to $\{F_1,\dots,F_{N-2},\,F_{N-1}\}$, we have that
\begin{equation}\label{cisn}
c_i^{(N)} = c_i^{(N-1)} + \ds\frac{1}{2} \alpha_i  \alpha_{N-1} \rho_{i,N-1} \ds\frac{\beta}{F_{N-1}^{1-\beta}} F_i^{\beta} V^2,\quad i=1,\dots,N-2.
\end{equation}
Since usually $F_i^{max}=F^{max},\,i=1,\dots,N-1$ and $V^{max}>>F^{max}$,
we can see that the advection increases with $N$ whereas the diffusion coefficients do not.

A way to avoid the boundary layers that this advection-dominance can produce with conditions \eqref{bound1_0} when $0<\beta \le 1$, is to consider homogeneous Neumann conditions instead \cite{hundsdorfer03nso}, i.e.,
\begin{equation}\label{bound1}
\begin{array}{l}
u(t,F_1,\dots,F_{N-1},F_N)=g(T,F_1,\dots,F_{N-1}), \quad \hbox{if } F_j=0 ,\,1\le j\le N , \\[0.6pc]
 \ds\frac{\partial u}{\partial F_j}(t,F_1,\dots,F_{N-1},F_N)=0,  \quad \hbox{if }\, F_j=F_j^{max},\,1\le j\le N   -1,\, \hbox{or }\,  F_N=V^{max}.
\end{array}
\end{equation}

On one hand, these homogeneous Neumann boundary conditions are appropriate from the financial point of view. Actually, if we analyse the behaviour of the payoff of the swaption $T_1\times(T_N-T_1)$, which is given by
$$g(T,F_1,F_2,\dots,F_{N-1})=\max \left\{ \sum_{i=1}^{N-1}\frac{\tau_i(F_i-K)}{(1+\tau_1 F_1)\cdots(1+\tau_i F_i)},0\right\},$$
when some forward rate $F_j\rightarrow \infty$, we can consider the approximation
$$
 g(T,F_1,\dots,F_{N-1},V) \approx f(F_1,\dots,F_{N-1}):=\ds\sum_{i=1}^{N-1} \ds\frac{\tau_i(F_i-K)}{P_i},
$$
where $P_i:=\ds\prod_{l=1}^i (1+\tau_l F_l),\quad 1\le i\le N-1$, and prove that
\begin{equation}\label{limit1}
\ds\frac{\partial f}{\partial F_j}(T,F_1,\dots,F_{N-1},V) \longrightarrow 0 \quad \hbox{when } F_j \rightarrow \infty,\quad j=1,\dots,N-1.
\end{equation}
In order to prove this, it is enough to write the partial derivatives of $f$ as
\begin{align*}
 & \ds\frac{\partial f}{\partial F_j}(F_1,\dots,F_{N-1})= \\
 & \qquad \ds\frac{\tau_j}{(1+\tau_1F_1) \cdots (1+\tau_{j-1}F_{j-1})(1+\tau_jF_j)^2} \left( (1+\tau_jK) - \ds\sum_{i=j+1}^{N-1} \ds\frac{\tau_i(F_i-K)}{\tilde{P}_{ij}} \right),
\end{align*}
where $ \tilde{P}_{ij}=P_i/P_j=(1+\tau_{j+1}F_{j+1})\cdots (1+\tau_iF_i),\,i\ge j+1$.
As in \cite{lopezvazquez18}, an homogeneous Neumann boundary condition is imposed at $F_N=V^{max}$ because the price of the derivative becomes independent of $V$ when $V$ approaches to infinity.

On the other hand, the choice \eqref{bound1} becomes better when some discretization of the ``spatial'' variables $F_j$ is used to approximate the solution of a PDE \eqref{adr} on a uniform spatial grid, as will be illustrated later on.


\section{Space discretization with finite differences}

Following the ideas given in \cite{lopezvazquez18}, firstly a space discretization is performed on a uniform spatial grid on $\Omega=[0,F_1^{max}]\times [0,F_2^{max}]\times \dots \times [0,F_{N-1}^{max}]\times [0,V^{max}]$. However, we propose here a different time discretization by using the AMFR-W-methods introduced in \cite{amfrw18}.

For the space discretization we consider $N$ integers  $(M_1,\dots,M_{N-1},M_N)$ to define on $\Omega$ the spatial grid with $M_i+1$ equally spaced points (denoting $F_N^{max}=V^{max}$) at the $F_i-$direction
$F_{i,j_i}=j_i h_i$, $0\le j_i\le M_i$, $h_i= \ds\frac{F_i^{max}}{M_i}$, $1\le i\le N$,
and discretize the derivatives in \eqref{mmNd} with second order central finite differences at each spatial node $(F_{1,j_1},\dots,F_{N-1,j_{N-1}},F_{N,j_{N}})$. More precisely, the  MoL  approximates the solution at each spatial point $u(t,F_{1,j_1},\dots,F_{N-1,j_{N-1}},F_{N,j_{N}})\approx U_{j_1,\dots,j_{N-1},j_N}(t)$, where the values $U_{j_1,\dots,j_{N-1},j_N}(t)$ need to satisfy the semi-discretized ODE system
\begin{equation}\label{mmnd_0}
\begin{array}{rl}
\ds\frac{d}{dt}U_{j_1,\dots,j_{N-1},j_N}(t)=&
 \ds\sum_{i=1}^{N}
d_i (F_{i,j_i},F_{N,j_N})  \,\Delta^{(i)}_{j_1,\dots,j_{N-1},j_N} \\[0.7pc]
&+\ds\sum_{i=1}^{N-1} \ds\sum_{k=i+1}^{N} m_{ik}  (F_{i,j_i},F_{k,j_k},F_{N,j_N}) \, \Delta^{(ik)}_{j_1,\dots,j_{N-1},j_N} \\[0.7pc]
&+ \ds\sum_{i=2}^{N-1} a_i (F_{2,j_2},\dots,F_{i,j_i},F_{N,j_N})  \, \nabla^{(i)}_{j_1,\dots,j_{N-1},j_N},
\end{array}
\end{equation}
with $\Delta^{(i)},\,\Delta^{(ik)},\nabla^{(i)}$ representing the approximations with central differences for the derivatives $(\partial^2 u/\partial F_i^2)$, $(\partial^2 u/\partial F_i\partial F_k)$, $(\partial u/\partial F_i)$, respectively.

One important drawback of finite differences comes from the complexity of dealing with all these approximations when the spatial dimension $N$ is large. For the efficient manipulation of these differences, we propose to use $N-$dimensional multi-indices $\mathbf{j}=(j_1,\dots,j_N)$, and the following Lemma \ref{lemmabij} that is proved in the Appendix.

\begin{lemma}\label{lemmabij}
Given $N$ pairs of integers $m_i \le M_i$ for $i=1,\dots,N$, let us define the set
${\cal I}_N=\{ \mathbf{j}=(j_1,\dots,j_N) \,\,|\,\, m_i\le j_i \le M_i,\,\,\forall i=1,\dots,N \}$.
Moreover, for $M_T=\prod_{k=1}^N (M_k-m_k+1)$, define the map
$\vartheta: {\cal I}_N \longrightarrow \{m_1,m_1+1,\dots,M_T+m_1-1\}$,
\begin{equation}\label{bij0}
 \vartheta (\mathbf{j})= J=j_1+\ds\sum_{l=2}^N \left( (j_l-m_l) \prod_{r=1}^{l-1} (M_r-m_r+1) \right).
\end{equation}
Then, the map $\vartheta$ is bijective.
\end{lemma}
Besides the proof of Lemma \ref{lemmabij}, in the Appendix it is also included a practical way for computing the inverse map $\vartheta^{-1} (J)=\mathbf{j}$ (see \eqref{invj} in Lemma \ref{algor1}).


 Note that due to the Dirichlet boundary conditions in \eqref{bound1}, the values for $j_i=0$ for $1\le i\le N$  are given by the derivative payoff. As a consequence, the ODE system \eqref{mmnd_0} is applied only when $j_i=1,\dots,M_i , \,1\le i\le N$,  so
it has dimension $L=M_1 \cdots M_{N-1} M_N$.
%
%
Then, we separate the multi-indices that correspond with finite differences nodes on the lower boundaries from the rest of them, so that two different bijections of type \eqref{bij0} are considered:
$$\begin{array}{l}
\vartheta_0: {\cal I}_N^{(0)}=\{ \mathbf{j}=(j_1,\dots,j_N) \,\,|\,\, 0\le j_i \le M_i,\,\,\forall i=1,\dots,N \} \longrightarrow \{0,1,\dots,M-1\},\\[0.5pc]
\vartheta_1: {\cal I}_N^{(1)}=\{ \mathbf{k}=(k_1,\dots,k_N) \,\,|\,\, 1\le k_i \le M_i,\,\,\forall i=1,\dots,N \} \longrightarrow \{1,2,\dots,L\},\\[0.5pc]
 \end{array}$$
where $M=\prod_{k=1}^N (M_k+1)$, $\,{\cal I}_N^{(1)}\subset {\cal I}_N^{(0)}$, and
$$\begin{array}{l}
\vartheta_0 (\mathbf{j})= J=j_1+\ds\sum_{l=2}^N \left( j_l \prod_{r=1}^{l-1} (M_r+1) \right),\quad \mathbf{j} \in {\cal I}_N^{(0)},  \\[0.5pc]
\vartheta_1 (\mathbf{k})= K=k_1+\ds\sum_{l=2}^N \left( (k_l-1) \prod_{r=1}^{l-1} M_r \right),\quad \mathbf{k} \in {\cal I}_N^{(1)}.
\end{array}$$
The set of integers $\{0,1,\dots,M-1\}$ can be obtained as the union of the two disjoint sets
$\Inn= \vartheta_0 \left({\cal I}_N^{(1)}\right)=\vartheta_0 \left(\vartheta_1^{-1}\left( \{1,\dots,L\} \right)\right),$ $\Out=\{0,1,\dots,M-1\}-\Inn$,
and we consider the vector $Y(t)=\left(Y_J(t)\right)_{J=0}^{M-1}$, where for each $J=0,\dots,M-1$, $(j_1,\dots,j_{N-1},j_N)=\vartheta_0^{-1}(J)$,  
\begin{equation}\label{dis0}
Y_J(t)= \left\{ \begin{array}{ll}  U_{j_1,\dots,j_{N-1},j_N}(t), & \hbox{if } J \in \Inn, \\[0.5pc]
g(T,F_{1,j_1},\dots,F_{N-1,j_{N-1}}), & \hbox{if } J \in \Out.
\end{array}\right.
\end{equation}
%
%
Note that $\vartheta_0(\mathbf{j}-\mathbf{e}_i)\in \{0,1,\dots,M-1\}$, for all
$J\in \Inn$ with $\mathbf{j}=\vartheta_0^{-1}(J)\in {\cal I}_N^{(1)}$. Moreover, $\vartheta_0(\mathbf{j}-\mathbf{e}_i)=J-E_i$ where
$E_1=1$, $E_i=\prod_{r=1}^{i-1}(M_r+1),\,\,\,i\ge 2$.
On the other hand, taking into account boundary conditions \eqref{bound1}, for all $J\in \Inn$ we have that $\mathbf{j}=\vartheta_0^{-1}(J)\in {\cal I}_N^{(1)}$ and $\vartheta_0(\mathbf{j}+\mathbf{e}_i)\in \{0,1,\dots,M-1\}$ except when $j_i=M_i$, since  $\mathbf{j}+\mathbf{e}_i=(\dots,M_i+1,\dots)\notin {\cal I}_N^{(0)}$.
So, if $J\in \Inn$, we have $J+ E_i=\vartheta_0(\mathbf{j}+\mathbf{e}_i)$, for all $i=1,\dots,N$. However, for the case $j_i=M_i$ we take the virtual value $Y_{J+E_i}=Y_{J-E_i}$ due to the homogeneous Neumann conditions \eqref{bound1} at these boundaries.
%
Therefore, taking into account \eqref{dis0}, the ODE system \eqref{mmnd_0} is shown in Sketch 1.
 \begin{sketch}{ODE system \eqref{mmnd_0}}
\begin{algorithmic}
\FOR{$J\in \Out$} 
\STATE {$Y_J'(t)=0$} 
\ENDFOR
\FOR{$J \in \Inn$} 
 \STATE {$\mathbf{j}= \vartheta_0^{-1}(J)=(j_1,\dots,j_N)$} 
   \FOR{$i=1,\dots,N$} 
  \STATE {$\left(d_i\right)_J=d_i (F_{i,j_i},F_{N,j_N})$; $\left(a_i\right)_J=a_i (F_{2,j_2},\dots,F_{i,j_i},F_{N,j_N})$}
  \STATE {$\Delta_J^{(i)} = \left\{\begin{array}{ll}
\ds\frac{Y_{J+E_i}-2Y_J+Y_{J-E_i}}{h_i^2}, & \hbox{if } j_i\ne M_i \\[0.7pc]
\ds\frac{2Y_{J-E_i}-2Y_J}{h_i^2},& \hbox{if }  j_i= M_i
\end{array}\right.$} 
\STATE {$\nabla_J^{(i)}= \left\{\begin{array}{ll}
\ds\frac{Y_{J+E_i}-Y_{J-E_i}}{2h_i}, & \hbox{if }  j_i\ne M_i\\
0, &  \hbox{otherwise}
\end{array}\right.$}  
 \ENDFOR
 \FOR{$1\le i\le N-1,\,\, i+1\le k\le N$} 
 \STATE {$\left(m_{ik}\right)_J= m_{ik}  (F_{i,j_i},F_{k,j_k},F_{N,j_N})$}
\STATE {$\Delta_J^{(ik)} =\left\{\begin{array}{ll}
 \ds\frac{Y_{J+E_i+E_k}+Y_{J-E_i-E_k}-Y_{J+E_i-E_k}-Y_{J-E_i+E_k}}{4h_i h_k}, &
\hbox{if }  j_i\ne M_i  \\ & \hbox{ and } j_k\ne M_k  \\[0.5pc]
0, &  \hbox{otherwise}
\end{array}\right.$} 

 \ENDFOR 
 
 \STATE {$Y_J'(t)= \ds\sum_{i=1}^{N} \left(d_i\right)_J  \, \Delta_J^{(i)}
+\ds\sum_{i=1}^{N-1} \ds\sum_{k=i+1}^{N} \left(m_{ik}\right)_J   \, \Delta_J^{(ik)}
+ \ds\sum_{i=2}^{N-1} \left(a_i\right)_J   \, \nabla_J^{(i)}$} 
\ENDFOR
\end{algorithmic}
 \end{sketch}

Thus, we obtain the semi-discretized {\it autonomous initial value problem (IVP)}
\begin{equation}\label{odeNd}
Y'= {\cal F}(Y),\qquad Y(0)=Y_0,\qquad t\in [0,T], \qquad {\cal F}(Y)={\cal F}_0(Y)+ \ds\sum_{i=1}^N{\cal F}_i(Y).
\end{equation}
In this splitting, for each $i=1,\dots,N$, the term ${\cal F}_i(Y)$ contains the second order differences
in the $F_{i}-$direction. That is, for all $J=0,1,\dots,M-1$, we have
\begin{equation}\label{ffis}
\left({\cal F}_i(Y)\right)_J= \left\{\begin{array}{ll}
\left(d_i\right)_J  \, \Delta_J^{(i)} , & \hbox{if } J\in \Inn, \,\, (j_1,\dots,j_N)=\vartheta_0^{-1}(J), \\[1pc] 0, & \hbox{if } J\in \Out,
\end{array}\right.
\end{equation}
%
while the term ${\cal F}_0(Y)$ gathers the discretization corresponding to the remain terms.

Clearly, $\forall i=0,1,\dots,N$, ${\cal F}_i(Y)={\cal A}_i Y$, where
for all $J\in \Out$, the $J$-th row of the constant matrix ${\cal A}_i$ is null. Besides, when $i\ge 1$ and $J\in \Inn$, the $J$-th row of ${\cal A}_i$ has at most three non-zero elements, located at the columns $J-E_i,\,J$ and $J+E_i$.
Therefore, the differential system $Y'= {\cal F}(Y)$ in \eqref{odeNd} involves actually $L$ unknowns since the corresponding equation for each $J\in \Out$ is null.
Once the semi-discretized linear IVP \eqref{odeNd} has been posed, the AMFR-W-methods given in \cite{amfrw18} can be applied for its time integration.


\section{Time discretization}

\subsection{AMFR-W-methods}
AMFR-W-methods belong to the class of W-meth\-ods \cite{steihaug79aat} for the time integration of IVPs of type \eqref{odeNd}. Thus, from an approximation $Y_n$ of the solution $Y(t)$ at $t=t_n$ and a step size ${\Delta t}>0$, an $s-$stage W-method gives the approximation
$Y_{n+1}$ at $t_{n+1}=t_n+{\Delta t} $ by
\begin{equation} \label{row-s} \begin{array}{rcl}
(I-\theta {\Delta t} W)K_r&=&\displaystyle {\Delta t} {\cal F} \Bigl( Y_n+\sum_{j=1}^{r-1} a_{rj} K_j\Bigr) +
\sum_{j=1}^{r-1} q_{rj}K_j, \quad r=1,2,\ldots,s,\\
Y_{n+1}&=&\displaystyle Y_n + \sum_{r=1}^s b_r K_r. \end{array}
\end{equation}
Each W-method is characterized by its coefficients $(A,Q,b,\theta)$, where
$A=(a_{rj})_{j<r}$, $Q=(q_{rj})_{j<r}$ and $b=(b_r)$, and by the arbitrary matrix $W$.
This kind of methods can also be understood as a generalization of Rosenbrock methods, which are obtained when
$W={\cal F}'(Y_n)$.
In order to get W-methods of high order, $W$ must be some rough approximation of ${\cal F}'(Y_n)$, and methods of order 3 and higher can be found in the literature under the assumption (see, for instance, \cite{gerisch02osa,rang05nrw,lang13wmi,gonzalez17wmt})
\begin{equation}\label{Wmfp}
W- {\cal F}'(Y_n)={\cal{O}}({\Delta t}),\qquad {\Delta t} \to 0.
\end{equation}

When ${\cal F}'(Y)$ admits a directional splitting of type \eqref{odeNd}, i.e. $ {\cal F}'(Y_n)=  {\cal F}_0'(Y_n) + \sum_{i=1}^N {\cal A}_i$, where the matrices ${\cal A}_i$ have simple structures, in \cite{amfrw18} the authors propose
\begin{equation}\label{Wmat1}
\Bigl(
 I - \theta {\Delta t} W\Bigr)^{-1} = \prod_{i=N}^1 \Bigl(
 I - \nu {\Delta t} {\cal A}_i\Bigr)^{-1}
\Bigl( 2I -\bigl(I - \theta{\Delta t} {\cal F}'(Y_n)\bigr)
 \prod_{i=N}^1 \Bigl(
 I -  \nu{\Delta t}  {\cal A}_i\Bigr)^{-1}\Bigr),
\end{equation}
thus combining the Approximate Matrix Factorization (AMF) technique from \cite{vanderhouwen01aff,hundsdorfer03nso} for the matrix $I-\theta{\Delta t} W$ with a refinement to the solution of the linear systems \cite[Section 3]{GHPjcam2014}.
With this selection of the matrix $W$, the condition (\ref{Wmfp}) is fulfilled and the introduction of a new parameter $\nu$ allows to improve the stability of the method.

Expanding these formulas on the semi-discretized IVP \eqref{odeNd}, each stage of the resulting AMFR-W-method \eqref{row-s}-\eqref{Wmat1} is computed as
\begin{equation}\label{amfrw}
\begin{array}{rcl}
K_r^{(0)}&=&{\Delta t} {\cal F}(Y_n+  \sum_{j=1}^{r-1} a_{rj} {K}_j) +
\sum_{j=1}^{r-1} q_{rj}{K}_j, \\[2mm]
(I-\nu {\Delta t}  {\cal A}_{i})K_r^{(i)}&=& K_r^{(i-1)}, \quad (i=1,\ldots,N)\\[2mm]
\hat{K}_r^{(0)} &=& 2 K_r^{(0)} - (I-\theta {\Delta t} {\cal F}'(Y_n)) K_r^{(N)},\\[2mm]
(I-\nu {\Delta t}   {\cal A}_i)\hat{K}_r^{(i)}&=& \hat{K}_r^{(i-1)}, \quad (i=1,\ldots,N)\\[2mm]
{K}_r&=&\hat{K}_r^{(N)}.
\end{array}
\end{equation}

In \cite{amfrw18}, different choices for the coefficients of these methods are tested. In this article we have used the proposed {\bf AMFR-W2}, that is a 2-stage AMFR-W-method, with coefficients \cite[p.\,400]{hundsdorfer03nso}
\begin{equation}\label{coef_W2}
a_{21}=2/3,\quad q_{21}=-4/3, \quad  b_1=5/4,\quad b_2=3/4.
\end{equation}
 Since  \eqref{Wmfp} is fulfilled, this method is of
order $3$ for $\theta=(3+\sqrt 3 )/6$ (in ODE sense). In \cite{amfrw18}, the authors also proved that the choice of the parameter $\nu$ depends on the number $N$ of terms in the splitting \eqref{odeNd} to get unconditional stability. More precisely, they apply this method on the parabolic test problem given in \cite{gonzalez17pde} and guarantee that this method is unconditionally stable on linear constant coefficients PDEs with mixed derivatives of dimension $N$ if $\nu \ge N \kappa_N \,\theta$ with the values of $\kappa_N$ given in \cite[Table 2]{amfrw18}, when both periodic and homogeneous Dirichlet boundary conditions.  In \cite{gonzalezHeston}, the authors solve the two-dimensional PDE for the well-known Heston model in options pricing. For this purpose, a hyperbolic change of variables is previously applied to the PDE, thus allowing the use of full non-uniform spatial meshes. However, we do not apply this change of variables since we use sparse grids to approximate efficiently the solution of \eqref{mmNd_no} for higher spatial dimensions.

Obviously, if the solution of linear systems in \eqref{amfrw} turns out too expensive from the computational point of view, the applicability of these schemes remains very limited. However, in the case of the PDE problem here addressed or similar multi-dimensional linear problems, due to the simple structure of the matrices ${\cal A}_i$, each linear system of type $(I-\nu {\Delta t}  {\cal A}_{i})K= G$ can be solved by using $\tilde{L}_i=\prod_{ k\ne i}^N M_k$ tridiagonal  linear systems of dimension $M_i$. In order to make easier the reading of this article, the details of this computation are included in Algorithm \ref{alg:solveLinearSystem} in the Appendix.

Moreover, another interesting advantage of the AMFR-W-methods \eqref{amfrw} when applied to the autonomous linear problem  \eqref{odeNd} with ${\cal F}(Y)={\cal A} Y$, ${\cal A}={\cal A}_0+\sum_{i=1}^N {\cal A}_i$, comes from the fact that the matrix-vector product ${\cal F}'(Y_n) K_r^{(N)}$ is simply an extra evaluation of the derivative function ${\cal F}(K_r^{(N)})$ and the explicit computation of the matrix ${\cal A}_0$ is not actually necessary.

\subsection{$\theta-$method + Gauss-Seidel as an $W$-method}

In \cite{lopezvazquez18}, the authors applied a direct (backwards in time) time-space discretization with finite differences, that can be also interpreted as a W-method \eqref{row-s} with fixed time step-size, when a fixed number of iterations of the Gauss-Seidel iterative scheme is used to solve the involved linear systems.

More precisely, for a fixed time step-size ${\Delta t}>0$, the well-known $\theta-$method applied to \eqref{odeNd}, with
$\theta \in [0,1]$, gives the approximations $W_n \approx Y(t_n),\, t_n=n {\Delta t},$ $n=0,1,\dots,M$, by using the formula
\begin{equation}\label{thetamet}
W_{n+1}= W_n + (1-\theta) {\Delta t}  {\cal F}(W_n)+ \theta  {\Delta t} {\cal F}(W_{n+1}).
\end{equation}
Therefore, when $\theta\ne 0$ and linear problems ${\cal F}(Y)={\cal A}Y$ as \eqref{odeNd} are considered, it is necessary to solve the linear system
$(I-\theta{\Delta t} {\cal A}) W_{n+1}= \beta_n$, at each time step, with $\beta_n=(I-\theta{\Delta t} {\cal A}) W_{n}+{\Delta t} {\cal F}(W_{n})$. In \cite{lopezvazquez18}, the Gauss-Seidel iterative linear systems solver is performed
 until getting an error below a prescribed tolerance. Note that Gauss-Seidel method splits the coefficient matrix ${\cal A}=P+R$, where $P$ is the triangular matrix whose entries are the lower-triangular part of ${\cal A}$ and its diagonal elements, while $R$ stores its strictly upper-triangular part. By using this splitting, from a starting value $W_{n+1}^{(0)}$, this method computes iteratively approximations $W_{n+1}^{(r)}\approx W_{n+1}$ by solving only triangular systems
\begin{equation}\label{GS}
(I-\theta{\Delta t} P) W_{n+1}^{(r)}=\theta{\Delta t} R W_{n+1}^{(r-1)} +\beta_n, \qquad  r=1,2,\dots
\end{equation}
After some algebraic manipulations and taking as a natural choice for the starting value $W_{n+1}^{(0)}=W_n$, this iterations can be written as $ W_{n+1}^{(r)}= W_n+ \ds\sum_{j=1}^r \hat{K}_j$, $r=1,2,\dots,$ where the vectors $\hat{K}_j$ are sequentially computed by
$$(I-\theta{\Delta t} P) \hat{K}_r = - (I-\theta{\Delta t} {\cal A}) \left(\sum_{j=1}^{r-1} \hat{K}_j\right) +{\Delta t} {\cal F}(W_{n}), \qquad r=1,2,\dots$$
or, equivalently
$$(I-\theta{\Delta t} P) \hat{K}_r = {\Delta t} {\cal A} \left(W_n+\sum_{j=1}^{r-1} \theta \hat{K}_j\right)+ \sum_{j=1}^{r-1} (-1) \hat{K}_j \qquad r=1,2,\dots$$
If we compare this last formula with \eqref{row-s}, clearly if this combination of $\theta$-method + Gauss-Seidel iteration \eqref{thetamet}-\eqref{GS} is performed with a fixed number $s$ of iterations, then the method can be included in the class of W-methods \eqref{row-s} with coefficients $W=P$, $a_{rj}=\theta, \,\forall j<r$, $q_{rj}=-1,\,\forall j<r$, $b_r=1,\, r=1,\dots,s$.
Therefore, if we understand the discretization proposed in \cite{lopezvazquez18} in the W-methods framework, we can apply the order conditions given in \cite[p.115]{hairer96sod} or \cite[Sec.2.1]{gonzalez16afo}  (in a similar notation as here), and we obtain that it achieves order 2 (for $s\ge 2$) in time only when $\theta=1/2$ (Crank-Nicolson scheme), what is in agreement with the results obtained in the aforementioned article \cite{lopezvazquez18}.

An advantage of expressing the scheme in \cite{lopezvazquez18} as a W-method is that it makes easier to compare its computational cost per time step with that of the here proposed {\bf AMFR-W2} method \eqref{amfrw}-\eqref{coef_W2}. The scheme \eqref{thetamet}-\eqref{GS} with $\theta=1/2$ and $s$ Gauss-Seidel iterations needs to compute one evaluation of the derivative function ${\cal F}$ and $s$ triangular linear systems of dimension $L=M_1\cdots M_N$. On the other hand, {\bf AMFR-W2} needs to evaluate four times the derivative function  (as ${\cal F}'(Y_n) K_r^{(N)}={\cal F}(K_r^{(N)}) )$ and to solve $(2\,\hat{L}_i)$ tridiagonal linear systems of dimension $M_i$, per each $i=1,\dots,N$. Since the triangular systems of dimension $L$ need approximately ${\cal O}(L^2)$ operations, the tridiagonal ones of dimension $M_i$ cost ${\cal O}(M_i)$ operations and each evaluation of derivative function involves ${\cal O}(L^2)$ operations, we have ${\cal O}((s+1)L^2)$ operations for \eqref{thetamet}-\eqref{GS} and
${\cal O}(4L^2+2NL)$ operations for {\bf AMFR-W2} \eqref{amfrw}-\eqref{coef_W2}.


\section{Sparse grids in space} \label{sec:standardSparseGrid}
Solving PDE problems as \eqref{mmNd_no} on a full tensor product based grid with $p^N$ grid points, with $p$ being the number of grid points in each coordinate direction, can become a highly involved computational task, even prohibitive. As the number of underlying forward rates increases, clearly the dimension of the multi-dimensional pricing PDE \eqref{mmNd_no} increases as well, so the computational cost of solving the fully discretized problem grows exponentially. Thus, the discretization using this so-called full grid also consumes too much memory. This drawback is referred as the curse of dimensionality. For example, pricing a swaption over five forward rates ruled by the same stochastic volatility, by means of a full grid with $128$ points per coordinate gives rise to more than four trillion points. The storage of such a grid using double precision floating point format will need more than $32$ thousand gigabytes of memory.

Because of the curse of dimensionality, traditional full grid methods, like finite differences, finite elements or finite volumes, are not able to price derivatives with high dimensional underlying processes, even in the most powerful supercomputers available nowadays. This limitation can be partially overcome by using a family of techniques known as sparse grid methods (see \cite{bungartzGriebel2004}, for example). Sparse grids are useful numerical methods for solving high-dimensional PDEs because they are based on a relatively small number of grid points but also maintain a satisfactory accuracy. More precisely, let $d$ denote the underlying problem's dimensionality and $p$ the number of grid points in one coordinate direction at the boundary. On the one hand, regarding the considered number of degrees of freedom, full grid methods use $O(p^d)$ grid points, while sparse grid discretizations only employ $O\left(p (\log_2 p)^{d-1}\right)$ grid points. On the other hand, concerning accuracy, conventional methods converge at a rate of $O(p^{-2})$ when making use of second order schemes, whereas sparse grid methods converge at the only slightly deteriorated rate of $O(p^{-2} \left(\log_2 p \right)^{d-1})$. In \cite{bungartzGriebel2004}, Bungartz and Griebel present an excellent survey of the fundamentals and the applications of sparse grids, with a focus on the solution of PDEs. Sparse grid were introduced in the early 1990s for the solution of PDEs by Zenger \cite{zenger91} and Griebel \cite{griebel91}.

\subsection{Standard sparse grid combination technique}

Discretizations on spar\-se grids require hierarchical data structures. Therefore, specially designed PDE solvers are required, and their implementations become more and more complicated as the dimension of the problem increases \cite{sAchatz, aZeiser}. An efficient way to avoid intricate sparse grid implementations is given by the sparse grid combination technique, originally proposed by Griebel, Schneider and Zenger \cite{griebelSchneiderZenger91}. Basically, the combination technique solves the PDE on several independent and conventional Cartesian smaller-sized grids. Then, the solution in the sparse grid space is approximated by a suitable linear combination of these partial solutions on the coarser grid. This solution retains the advertised convergence rate of sparse grid methods if certain error expansions for the component approximations exist \cite{bgrz94,bgrz96,aa04}. Note the rigorous analysis of finite differences schemes for the sparse grid combination technique in \cite{reisinger13}. Further advantages of the technique are the possibility to solve the problem on each of the constituent grids using standard full grid solvers and the inherent parallelism of the method \cite{griebel92}.

Let us introduce formally the sparse grid combination technique. We fix a multi-index $\pmb{l} = (l_1,l_2,\ldots,l_d)\in\mathbb{N}^d_0$ and define its $L_1$-norm as $|\pmb{l}|_1 = \sum_{i=1}^d l_i.$ In the $d$-dimensional orthohedron $[0,c_1]\times[0,c_2]\times \ldots \times [0,c_d]$  ($c_i\in \mathbb{R}_{>0}$), we denote by $\Omega_{\pmb{l}} = \Omega_{(l_1,\ldots,l_d)}$ an anisotropic\footnote{Mesh spacing differs in each coordinate direction.} although full grid having uniform mesh spacing $h_i = 2^{-l_i}c_i$ in each coordinate direction $i \in \{1,\ldots,d\}$. Let $u_{\pmb{l}}$ be the conventional finite difference solution to the PDE on grid $\Omega_{\pmb{l}}$, extended to $[0,c_1]\times[0,c_2]\times \ldots \times [0,c_d]$ by interpolation. Then, the sparse grid combination solution $u_n^s$ over the sparse grid $\Omega_n^s$  with refinement level $n$ is given by the following linear combination
\begin{equation}
 \label{eq:standardCombTechnique}
  u^s_n = \displaystyle\sum_{q=0}^{d-1} (-1)^q \cdot {d-1 \choose q} \cdot \displaystyle \sum_{|\pmb{l}|_1 = n-q} u_{\pmb{l}}.
\end{equation}
Increasing the level $n$ should give a more accurate solution to the problem. The grid solutions $u_\mathbf{l}$ involved in the inner sum of \eqref{eq:standardCombTechnique} all have $\mathbf{l}=(l_1,\ldots,l_d)$ such that $l_1+\cdots+l_d = n-q$. The number of elements in each of these grids is $O(2^{n-q})$, regardless of the dimension, and the number of grid solutions in this inner sum is ${n-q+d-1 \choose d-1}$ and grows like $O((n-q)^{d-1})$.
Besides, $\left( \bigcup\limits_{0\leq|\mathbf{l}|_1< n} \Omega_\mathbf{l} \right) \subset \Omega^s_n$. Therefore, the dimension of the sparse grid space on level $n$ is $O(2^n n^{d-1}) = O(h^{-1} (\log_2 h)^{d-1})$ where $h=2^{-n}$ is the finest grid size. This value can be compared with the full grid space dimension which is $O(2^{nd}) = O(h^{-d})$.




The combination technique works due to the cancellation mechanism of the error terms in the involved grids. This cancellation principle is well known in extrapolation techniques. Indeed, all lower order error terms cancel out in the the combination formula \eqref{eq:standardCombTechnique}, see \cite{hendricksEhrhardtGunther2017} for deep details in dimension two. Thus, the combination technique is able to produce accurate results in reasonable time. Several generalizations of the standard combination technique formula \eqref{eq:standardCombTechnique} have been developed \cite{hegland2007}.

The combination technique algorithm is embarrassingly parallel since all component grid solutions can be computed in parallel. In general, for refinement level $n$ in $d$ dimensions there are $\sum_{q=0}^{d-1} {n-q+d-1 \choose d-1}$ component grids, which can be solved in parallel. In order to achieve optimal speed-ups one has to carefully deal with load imbalances, even in all those grids at the same refinement level, due to the anisotropic structure of the component grids.

The sparse grid combination technique was initially formulated for elliptic PDEs such as Laplace's and Poisson's equation. Later, it has also been applied to parabolic PDE, specially for option pricing problems in finance \cite{aa04,leentvaarOosterlee2006,ReisingerWittum2007, leentvaarOosterlee2008, reisinger13, chiarellaKang2013, beyna2013, salasVazquez2017_477,lopezvazquez18,hendricksHeuerEhrhardtGunther2017, duringHendricksMiles,hendricksEhrhardtGunther2016,hendricksEhrhardtGunther2017}. Here we focus on the implementation of the sparse grid combination technique for parabolic equations. More precisely, we just consider the case where the solution is only needed at the final time, which is frequently the case in finance and particularly in the problem we address.

In this setting, the natural approach is the following. First, we solve the parabolic equation on each of the full grids involved in the sparse grid combination technique formula \eqref{eq:standardCombTechnique} with a full grid method. Finally, combine these solutions only at the end. This method only requires interpolation from grid values at the final time, but not at intermediate time steps. 
If the numerical error due to the time discretization does not dominate the spatial error, we expect  a pointwise rate of convergence proportional to $O(p^{-2} \left(\log_2 p \right)^{d-1})$ for our AMFR-W scheme applied to problems with smooth enough initial and boundary data. It is important to notice that interpolation techniques are required in order to approximate the solution at points not belonging to the sparse grid. The most straightforward approach is to interpolate at those points over all full grids handled by the combination technique, and then add up these results with the appropriate combination technique weights. Note also that the interpolation technique has to preserve the order of the used discretization scheme, so that the convergence result remain valid for the entire domain. Otherwise, the convergence order only holds for grid nodes belonging to all sub-grids and therefore not affected by interpolation. A tensor based linear interpolation preserves the required order $2$ of accuracy for second order finite difference discretizations.


\subsection{Modified sparse grid combination technique}

By means of the previously described standard sparse grid combination technique it is impossible to approximate accurately a Neumann boundary condition for degenerated Cartesian grids having very few points in the corresponding coordinate direction. 
The approximations of the solution in these grids becomes very poor, thus decreasing the accuracy of the combination technique approximation.

In order to overcome this drawback, a mild modification of the standard sparse grid combination technique \eqref{eq:standardCombTechnique} can be developed, just by forcing a minimum number of discretization steps in all grids involved in the combination procedure. More precisely, all levels in all dimensions start from a small but non zero value $\psi$, so that the modified combination technique formula reads \begin{equation}
 \label{eq:modifiedCombTechnique}
  u^s_n = \displaystyle\sum_{q=0}^{d-1} (-1)^q \cdot {d-1 \choose q} \cdot \displaystyle \sum_{|\pmb{l}|_1 = n-q} u_{\psi\mathbf{1} + \pmb{l}},
\end{equation}
where $\psi \mathbf{1} + \pmb{l} = (\psi+ l_1,\ldots, \psi + l_d)$.  This modified sparse grid combination technique working over a modified sparse grid $\Omega^{s,\psi}_n$ produces more accurate approximations \cite{beyna2013} at the cost of increasing the consumed time and memory. 
Although the number of subproblems to be solved is exactly the same as before, the number of degrees of freedom associated to each subproblem increases. In fact, the number of grid points in the combined sparse grid increases from $O( 2^n n^{d-1})$ in the standard one to $O(2^{n+d\psi} n^{d-1})$ in the modified one. In this new setting $\psi$ should be kept small (specifically $\psi=1$ or $2$ in the present work), otherwise the new modified combination technique will suffer soon the curse of dimensionality.





\section{Numerical results}
In this Section we present the obtained numerical results when the previously described methodologies are applied. More precisely, we show and discuss the results obtained by using the AMFR-W method with full grid, standard and modified sparse grids combination techniques to conveniently cope with the proposed homogeneous Neumann boundary conditions in the particular case $\beta=1$.

For all products we will use the same hypothetical market data presented in Table \ref{table:marketData} where we consider the tenor structure $0=T_0 < 0.5 < 1.0 < \cdots < 4.5 < 5 = T_{10}$ in years, with constant periods $\tau = T_{i+1}-T_i = 0.5$.

\begin{table}[!htb]
\begin{center}
{\footnotesize
\begin{tabular}{|r|r|r|r||r|r|r|r|}
\hline
$i$ & $T_i$ & $F_i(0)$ & $\alpha_i$ & $i$ & $T_i$ & $F_i(0)$ & $\alpha_i$\\
\hline
\hline
 $0$ & $0$ & $0.0112$ & $0$ & $3$ & $1.5$ & $0.0126$ & $0.2221$\\ %
\hline
 $1$ & $0.5$ & $0.0118$ & $0.2366$ & $4$ & $2$ & $0.0130$ & $0.2068$\\ %
\hline
 $2$ & $1$ & $0.0122$ & $0.2145$ & $5$ & $2.5$ & $0.0135$ & $0.1932$ \\ %
\hline
\end{tabular}
}
\end{center}
\caption{Hypothetical market data (LIBOR rates and volatilities) used in pricing. Strike rate $K = 0.011$.}
\label{table:marketData}
\end{table}

The spatial domain is defined by $F^{max} = 0.04$ and $V^{max} = 3.5$, thus upper boundaries were settled between $3$ and $4$ times the point of interest at which we evaluate the pricing of the interest rate derivative.
In the cases where the analytical solution is not available, we first compute reference solutions using the proposed space and time discretizations over classical full grids. These solutions will serve to assess about the accuracy of the proposed sparse grids methods in space.
The designed algorithms were implemented using C++ (GNU C++ compiler 8.3.1) and double precision. Besides, all numerical experiments have been performed in a machine with 16GBytes of RAM and four multicore Intel Xeon CPUs E5-2620 v4 clocked at 2.10GHz, each one with eight cores. 

\subsection{Numerical results with full grids}

The first test to validate the proposed numerical methodologies consists of pricing a caplet without considering stochastic volatility, that is to say, under the classical LMM. This test is a sanity check, since the analytical pricing formula is known for caplets, the so-called Black-Scholes's formula for caplets(\cite[equation 1.26]{brigo_mercurio2007}). More precisely, we start pricing the caplet with maturity $T_1$ and payoff $\tau_1(F_1(T_1)-K)^+$ paid at time $T_2$, under the data of Table \ref{table:marketData}, with strike rate $K=0.011$. The present intrinsic value of the caplet is given by $P(0,T_2)\tau_1(F_1(T_1)-K)^+.$ The exact price of this product given by Black-Scholes' formula is $6.058877$ basis points (bps, $1 \mbox{ bp }=10^{-4}$).

In order to price this caplet using the here presented PDE approach it is convenient to consider the terminal probability measure associated with choosing the bond $P(0,T_2)$ as numeraire. Thus, the price of this product is given by the solution of the PDE \eqref{mmNd_no} (with $\sigma=0$) multiplied by $P(0,T_2) = \frac{1}{1+\tau F_0(0)}\frac{1}{1+\tau F_1(0)}$, i.e. $P(0,T_2) u(0,F_1,V)$. Once obtained the PDE solution on the last time slice, the price of the caplet is obtained through interpolation in space, by means of multilinear interpolation, thus maintaining order two in space. In Table \ref{table:fg_sigma0} full grid solutions are presented for levels from $6$ to $13$ in space and considering $4$, $8$, $16$ and $256$ time steps. The interpolation in space for the last time slice was done in $F_1=0.0118$ and $V=1$. The column labelled as Solution shows the PDE solution in bps, and the column for the error measures the absolute distance of the numerical solution to the exact one, in bps as well. The execution time was measured in seconds in all the experiments in this work. The ``grid points'' column displays the number of grid points employed in the full space meshes at each time discretization. Since the method is order three in time, few time steps could be considered in real pricing applications. For the space level $13$ (for the forward and the volatility) and when using $256$ time steps, the method was able to recover the exact solution up to the $8$-th decimal digit. Nevertheless, in this case the full grid method required almost three hours, all space meshes in all time slices with more than $67$ million points.

\begin{table}[!htb]
\begin{center}
{\tiny
\begin{tabular}{|r|rrr|rrr|r|}
\hline
& \multicolumn{3}{c|}{$4$ time steps} & \multicolumn{3}{c|}{$8$ time steps} &\\
 \cline{2-8}
Level & Solution & Error & Time & Solution & Error & Time & Grid points\\
\hline
\hline
\hline
 $6$ & $6.084214$ & $2.533666\times 10^{-2}$ & $0.01$ & $6.082708$ & $2.383026\times 10^{-2}$ & $0.02$ & $4225$ \\ %
\hline
 $7$ & $6.065957$ & $7.079711\times 10^{-3}$ & $0.05$ & $6.064274$ & $5.396418 \times 10^{-3}$ & $0.09$ & $16641$ \\ %
\hline
 $8$ & $6.063526$ & $4.648521 \times 10^{-3}$ & $0.15$ & $6.061980$ & $3.102002 \times 10^{-3}$ & $0.27$ & $66049$ \\ %
\hline
 $9$ & $6.060148$ & $1.270262\times 10^{-3}$ & $0.55$ & $6.058939$ & $6.133732\times 10^{-5}$ & $0.96$ & $263169$ \\ %
\hline
 $10$ & $6.060300$ & $1.422802 \times 10^{-3}$ & $2.14$ & $6.059077$ & $1.992885 \times 10^{-4}$ & $4.23$ & $1050625$  \\ %
\hline
 $11$ & $6.060237$ & $1.359825\times 10^{-3}$ & $9.39$ & $6.059021$ & $1.431228 \times 10^{-4}$ & $18.69$ & $4198401$ \\ %
\hline
 $12$ & $6.060236$ & $1.358165 \times 10^{-3}$ & $40.47$ & $6.059019$ & $1.416918 \times 10^{-4}$ & $80.45$ & $16785409$ \\ %
\hline
 $13$ & $6.060226$ & $1.348640 \times 10^{-3}$ & $165.64$ & $6.059011$ & $1.331767\times 10^{-4}$ & $328.15$ & $67125249$ \\ %
\hline
\hline
& \multicolumn{3}{c|}{$16$ time steps} & \multicolumn{3}{c|}{$256$ time steps} & \\
 \cline{2-8}
Level & Solution & Error & Time & Solution & Error & Time & Grid points\\
\hline
\hline
\hline
 $6$ & $6.082540$ & $2.366282\times 10^{-2}$ & $0.04$ & $6.082513$ & $2.363549 \times 10^{-2}$ & $0.47$ & $4225$ \\ %
\hline
 $7$ & $6.064109$ & $5.231082\times 10^{-3}$ & $0.16$ & $6.064079$ & $5.201370\times 10^{-3}$ & $1.89$ & $16641$ \\ %
\hline
 $8$ & $6.061870$ & $2.992441\times 10^{-3}$ & $0.49$ & $6.061840$ & $2.962119\times 10^{-3}$ & $7.28$ & $66049$ \\ %
\hline
 $9$ & $6.058832$ & $4.556637\times 10^{-5}$ & $1.97$ & $6.058802$ & $7.582048 \times 10^{-5}$ & $30.77$ & $263169$ \\ %
\hline
 $10$ & $6.058975$ & $9.701419 \times 10^{-5}$ & $8.33$ & $6.058944$ & $6.680796\times 10^{-5}$ & $131.75$ & $1050625$ \\ %
\hline
 $11$ & $6.058919$ & $4.157937\times 10^{-5}$ & $37.19$ & $6.058889$ & $1.139059 \times 10^{-5}$ & $593.13$ & $4198401$ \\ %
\hline
 $12$ & $6.058918$ & $4.041366 \times 10^{-5}$ & $163.01$ & $6.058888$ & $1.022896 \times 10^{-5}$ & $2558.92$ & $16785409$ \\ %
\hline
 $13$ & $6.058910$ & $3.191855\times10^{-5}$ & $655.84$ & $6.058879$ & $1.735259 \times 10^{-6}$ & $10415.57$ & $67125249$ \\ %
\hline
\end{tabular}

}
\end{center}
\caption{Full grid method, caplet with expiry $T_1$, $\sigma=0$. Prices and errors are shown in bps.}
\label{table:fg_sigma0}
\end{table}

Once we have checked the correct behaviour of the full grid method, which will be used in the sparse grid combination technique, we compute full grid reference solutions for financial products without exact prices. They will be used in order to assess on the correctness of the oncoming sparse grid combination technique implementation. Therefore, in order to minimize errors due to the time discretization, $256$ time steps will be chosen for the rest of the full grid tests in this section.
In Table \ref{table:fg_2d} the computed prices of the previous caplet under the stochastic volatility framework are shown.

Next, we deal with the pricing of $T_a \times (T_b - T_a)$ European swaptions. In Table \ref{table:fg_3d}, first the results for the $0.5\times 1$ swaption are given. Note that under this full grid framework it is not possible to price this product in reasonable computational times past refinement level $9$, due to the high number of involved spatial grid points. Then, the results for the $0.5 \times 1.5$ swaption are also shown. Once more, full grid pricing is only achievable on the lower grid levels.

\begin{table}[!htb]
\begin{center}
{\footnotesize
\begin{tabular}{|r|r|r|r||r|r|r|r|}
\hline
Level & Solution & Time & Grid points & Level & Solution & Time & Grid points \\
\hline
\hline
 $6$ & $6.050103$ & $0.48$ & $4225$ & $10$ & $6.023799$ & $135.29$ & $1050625$ \\ %
\hline
 $7$ & $6.029510$ & $1.93$ & $16641$ & $11$ & $6.023737$ & $597.82$ & $4198401$ \\ %
\hline
 $8$ & $6.026929$ & $7.51$ & $66049$ & $12$ & $6.023734$ & $2557.62$ & $16785409$ \\ %
\hline
 $9$ & $6.023665$ & $31.06$ & $263169$ & $13$ & $6.023725$ & $10505.02$ & $67125249$ \\ %
\hline
\end{tabular}
}
\end{center}
\caption{Full grid method, caplet with expiry $T_1$, $\sigma=0.3$, $\phi_1=0.4$, $256$ time steps.}
\label{table:fg_2d}
\end{table}

\begin{table}[!htb]
\begin{center}
{\footnotesize
\begin{tabular}{|r|r|r|r|r||r|r|r|r|}
\hline
\multirow{3}{*}{\rotatebox[origin=c]{90}{$0.5\times 1$}} & Level & Solution & Time & Grid points & Level & Solution & Time & Grid points\\
\cline{2-9}
 & $6$ & $13.002003$ & $71.33$ & $274625$ & $8$ & $12.981320$ & $5111.14$ & $16974593$\\ %
\cline{2-9}
 & $7$ & $12.984709$ & $590.97$ & $2146689$ & $9$ & $12.980459$ & $43325.53$ & $135005697$\\ %
\hline
\hline
\multirow{3}{*}{\rotatebox[origin=c]{90}{{\tiny$0.5 \times 1.5$}}} & Level & Solution & Time & Grid points & Level & Solution & Time & Grid points\\
\cline{2-9}
 & $3$ & $23.952705$ & $1.70$ &  $6561$ & $5$ & $21.577149$ & $474.93$ &  $1185921$ \\
 \cline{2-9}
 & $4$ & $21.765442$ & $28.45$ &  $83521$ &  $6$ & $21.486079$ & $8122.91$ &  $17850625$\\ %
\hline
\end{tabular}
}
\end{center}
\caption{Full grid method, $0.5\times 1$ and $0.5\times 1.5$ swaptions, $\sigma=0.3$, $\phi_i=0.4$, $i=1,2,3$, $256$ time steps.}
\label{table:fg_3d}
\end{table}

\subsection{Numerical results with the standard sparse grid combination technique}

In this section, by means of the standard sparse grid combination technique, we price not only the previous caplets and swaptions, but also swaptions involving more underlying forward interest rates, thus dealing with high dimensional setting. As usual, we are also interested in the values of these derivatives at the last time cut for the values of the forward rates depicted in Table \ref{table:marketData} and $V=1$, which define the spatial point where the value of the solution of the PDE is computed. In order to obtain the solution given by the sparse grid combination technique at this point, the numerical solution on each grid involved in the combination technique is interpolated at this point with multilinear interpolation. Next, all these values are introduced in the combination technique formula \eqref{eq:standardCombTechnique}, thus obtaining the price provided by the standard sparse grid combination technique.

Moreover, sparse grid combination techniques have been implemented to take advantage of shared memory parallel computers. The code was optimized and parallelized using OpenMP framework \cite{ref:openmp}, version 4.5. In order to deal with the previously mentioned load imbalances it is crucial to use a dynamic schedule to assign the involved full grids to threads. In this way, OpenMP assigns one grid to each thread. When the thread finishes, it will be assigned the next mesh that has not been assigned yet. The speedup of the parallelized version is almost equal to the number of available computing cores, in our case $32$. This optimal speedup is due to the fact that communication between processors only takes place at the final step in order to concentrate the solutions over all grids to a single scalar value.

In Table \ref{table:d2_1} we price the caplet with maturity $T_1$ under the framework without stochastic volatility, whose exact price was $6.058877$ basis points. The results in this Table are to be compared with those of Table \ref{table:fg_sigma0}. The accuracy of each solution is only slightly worse in this case, although the computing time is much lower due the much less number of involved grid points. For example, with the full grid approach, the solution using $256$ time steps and refinement level in space $13$ employed $10415.57$ seconds to achieve an error $1.735259\times 10^{-6}$ in basis points, while the standard sparse grid combination technique just needed $569.26$ seconds to attain almost the same accuracy, an error of $1.047159 \times 10^{-5}$. The reduction in the number of employed grid points is also shown in Table \ref{table:fg_sigma0}.

\begin{table}[!htb]
\begin{center}
{\scriptsize
\begin{tabular}{|r|rrr|rrr|r|}
\hline
& \multicolumn{3}{c|}{$4$ time steps} & \multicolumn{3}{c|}{$8$ time steps} & \\
 \cline{2-8}
Level & Solution & Error & Time & Solution & Error & Time & Grid points\\
\hline
\hline
 $6$ & $6.063081$ & $4.202988 \times 10^{-3}$ & $0.02$ & $6.053883$ & $4.994378 \times 10^{-3}$ & $0.02$ & $385$ \\
\hline
 $7$ & $6.120850$ & $6.197245 \times 10^{-2}$ & $0.02$ & $6.110998$ & $5.212012 \times 10^{-2}$ & $0.02$ & $833$ \\
\hline
 $8$ & $6.067165$ & $8.287799 \times 10^{-3}$ & $0.02$ & $6.059333$ & $4.557323 \times 10^{-4}$ & $0.03$ & $1793$ \\
\hline
 $9$ & $6.061071$ & $2.192917 \times 10^{-3}$ & $0.03$ & $6.056314$ & $2.563871 \times 10^{-3}$ & $0.05$ & $3841$ \\
\hline
 $10$ & $6.061681$ & $2.803747 \times 10^{-3}$ & $0.07$ & $6.059090$ & $2.119811\times 10^{-4}$ & $0.12$ & $8193$ \\
\hline
 $11$ & $6.059735$ & $8.569127 \times 10^{-4}$ & $0.20$ & $6.058133$ & $7.450762\times 10^{-4}$ & $0.35$ & $17409$ \\
\hline
 $12$ & $6.060156$ & $1.278608\times 10^{-3}$ & $0.58$ & $6.058787$ & $9.054229\times 10^{-5}$ & $1.21$ & $36865$ \\
\hline
 $13$ & $6.060211$ & $1.333075 \times 10^{-3}$ & $2.43$ & $6.058951$ & $7.311969\times 10^{-5}$ & $4.53$ & $77825$ \\
\hline
 $14$ & $6.060223$ & $1.345681 \times 10^{-3}$ & $9.00$ & $6.058998$ & $1.206373 \times 10^{-4}$ & $17.90$ & $163841$ \\
\hline
\hline
& \multicolumn{3}{c|}{$16$ time steps} & \multicolumn{3}{c|}{$256$ time steps} & \\
 \cline{2-8}
Level & Solution & Error & Time & Solution & Error & Time & Grid points\\
\hline
\hline
 $6$ & $6.052778$ & $6.099826 \times 10^{-3}$ & $0.02$ & $6.052738$ & $6.139440 \times 10^{-3}$ & $0.05$ & $385$ \\
\hline
 $7$ & $6.110144$ & $5.126590 \times 10^{-2}$ & $0.03$ & $6.110118$ & $5.124050 \times 10^{-2}$ & $0.09$ & $833$ \\
\hline
 $8$ & $6.059007$ & $1.291280 \times 10^{-4}$ & $0.04$ & $6.058984$ & $1.066738 \times 10^{-4}$ & $0.22$ & $1793$ \\
\hline
 $9$ & $6.056190$ & $2.687266 \times 10^{-3}$ & $0.07$ & $6.056164$ & $2.713620 \times 10^{-3}$ & $0.71$ & $3841$ \\
\hline
 $10$ & $6.059029$ & $1.514420 \times 10^{-4}$ & $0.20$ & $6.058998$ & $1.205596\times 10^{-4}$ & $2.31$ & $8193$ \\
\hline
 $11$ & $6.058034$ & $8.433494\times 10^{-4}$ & $0.64$ & $6.058003$ & $8.741697\times 10^{-4}$ & $8.47$ & $17409$ \\
\hline
 $12$ & $6.058698$ & $1.800588\times 10^{-4}$ & $2.33$ & $6.058667$ & $2.108407 \times 10^{-4}$ & $34.02$ & $36865$ \\
\hline
 $13$ & $6.058852$ & $2.520960\times 10^{-5}$ & $8.99$ & $6.058822$ & $5.565743\times 10^{-5}$ & $141.73$ & $77825$ \\
\hline
 $14$ & $6.058897$ & $1.981956\times 10^{-5}$ & $35.69$ & $6.058867$ & $1.047159 \times 10^{-5}$ & $569.26$ & $163841$ \\
\hline

\end{tabular}
}
\end{center}
\caption{Sparse grid combination technique, caplet with expiry $T_1$, $\sigma=0$.}
\label{table:d2_1}
\end{table}

Next, in Table \ref{table:d1_sg_sigma} the results for the previous caplet under the stochastic volatility framework are shown. These results are to be compared with those of Table \ref{table:fg_2d}.
Then, Tables \ref{table:sg_d3} and \ref{table:sg_d4} show the prices given by the standard  sparse grid combination technique for $0.5\times 1$ and $0.5\times 1.5$ swaptions under stochastic volatility. These results are to be compared with those of Table \ref{table:fg_3d}. Clearly, the standard sparse grid combination technique outperforms the full grid approach. Besides, the sparse method is able to cope with higher resolution levels, thus allowing to price successfully the $0.5\times 1.5$ swaption. Note that this was not possible with the full grid approach, see Table \ref{table:fg_3d}.

\begin{table}[!htb]
\begin{center}
{\footnotesize
\begin{tabular}{|r|r|r||r|r|r|}
\hline
Level & Solution & Time & Level & Solution & Time\\
\hline
\hline
 $6$ & $6.057668$ & $0.06$ & $11$ & $6.022082$ & $8.40$ \\ %
\hline
 $7$ & $6.095685$ & $0.10$ & $12$ & $6.023257$ & $33.96$ \\ %
\hline
 $8$ & $6.025848$ & $0.25$ & $13$ & $6.023595$ & $141.41$ \\ %
\hline
 $9$ & $6.018222$ & $0.70$ & $14$ & $6.023693$ & $569.07$ \\ %
\hline
 $10$ & $6.021834$ & $2.32$  \\ %
\cline{1-3}
\end{tabular}
}
\end{center}
\caption{Sparse grid combination technique, caplet with expiry $T_1$, $\sigma=0.3$, $\phi_1 = 0.4$, $256$ time steps.}
\label{table:d1_sg_sigma}
\end{table}

\begin{table}[!htb]
\begin{center}
{\footnotesize
\begin{tabular}{|r|r|r|r||r|r|r|r|}
\hline
Level & Solution & Time & Grid points & Level & Solution & Time & Grid points\\
\hline
\hline
 $8$ & $12.311172$ & $0.67$ & $8705$ & $12$ & $13.205324$ & $89.80$ & $219137$\\ %
\hline
 $9$ & $13.024747$ & $2.00$ & $19713$ & $13$ & $12.993536$ & $360.72$ & $483329$\\ %
\hline
 $10$ & $13.616333$ & $6.77$ & $44289$ & $14$ & $12.971783$ & $1399.46$ & $1060865$\\ %
\hline
 $11$ & $13.525821$ & $24.42$ & $98817$ & $15$ & $12.973900$ & $5755.68$ & $2318337$\\ %
 \hline
\end{tabular}
}
\end{center}
\caption{Sparse grid combination technique, $0.5 \times 1$ swaption, $\sigma=0.3$, $\phi_1 = \phi_2 = 0.4$, $256$ time steps.}
\label{table:sg_d3}
\end{table}

\begin{table}[!htb]
\begin{center}
{\footnotesize
\begin{tabular}{|r|rr|rr|r|}
\hline
 & \multicolumn{2}{c|}{$8$ time steps} & \multicolumn{2}{c|}{$256$ time steps} &  \\
 \cline{2-6}
  Level & Solution & Time & Solution & Time & Grid points\\
\hline
 $12$ & $21.935448$ & $8.50$ & $21.936574$ & $271.54$ &$1064961$ \\ %
\hline
 $13$ & $21.842901$ & $31.75$ & $21.844522$ & $998.51$ &$2439169$ \\ %
\hline
 $14$ & $21.609183$ & $119.97$ & $21.610055$ & $3809.31$ & $5550081$ \\ %
\hline
 $15$ & $21.707363$ & $461.28$ & $21.708001$ & $14866.96$ & $12554241$ \\ %
\hline
 $16$ & $21.519917$ & $1838.36$ & $21.516402$ & $59010.15$ & $28246017$ \\ %
\hline
 $17$ & $21.483062$ & $7315.95$ & $21.478930$ & $235912.03$ & $63242241$ \\ %
\hline
\end{tabular}
}
\end{center}
\caption{Sparse grid combination technique, $0.5\times 1.5$ swaption $\sigma=0.3$, $\phi_1 =\phi_2 =\phi_3= 0.4$.}
\label{table:sg_d4}
\end{table}

Finally, in Tables \ref{table:sg_d5} and \ref{table:sg_d6}, $0.5 \times 2$ and $0.5 \times 2.5$ swaptions are priced under stochastic volatility. The curse of dimensionality makes impossible to price these products with full grid approaches. In order to speedup the convergence of the sparse grid method, a useful technique is to consider a computational domain such that the point of interest is in the neighbourhood of the center of the domain. This strategy easily improves sparse grid results. In fact, in that region is where the sparse grid contains more points. Indeed, the central point belongs to all non degenerated grids involved in the standard sparse grid combination technique. The improvement in accuracy can be observed in Table \ref{table:sg_d6}, where the upper boundaries of the forward rates, $F^{max}$, were shrunken from $0.04$ to $0.02$.

\begin{table}[!htb]
\begin{center}
{\footnotesize
\begin{tabular}{|r|r|r|r|r|}
\hline
 & \multicolumn{2}{c|}{$4$ time steps} & \multicolumn{2}{c|}{$8$ time steps} \\
 \cline{2-5}
 Level & Solution & Time & Solution & Time \\
 \hline
 $14$ & $35.341806$ & $180.34$ & $35.346408$ & $360.57$ \\ %
\hline
 $15$ & $34.388334$ & $669.85$ & $34.425087$ & $1335.40$ \\ %
\hline
 $16$ & $32.115380$ & $2561.87$ & $32.122101$ & $5133.81$ \\ %
\hline
 $17$ & $30.639336$ & $10058.18$ & $30.641664$ & $20076.39$ \\ %
\hline
 $18$ & $30.881086$ & $40097.11$ & $30.918448$ & $80268.57$ \\ %
\hline
 $19$ & $30.822037$ & $239746.17$ & $30.797087$ & $479681.87$ \\ %
 \hline
\end{tabular}
}
\end{center}
\caption{Sparse grid combination technique, $0.5\times 2$ swaption, $\sigma=0.3$, $\phi_1 =\ldots =\phi_4= 0.4$.}
\label{table:sg_d5}
\end{table}

\begin{table}[!htb]
\begin{center}
{\footnotesize
\begin{tabular}{|r|r|r|r|r|r|r|}
\hline
\multicolumn{7}{|c|}{$F^{max} = 0.04$}\\
\hline
 & \multicolumn{2}{c|}{$2$ time steps} & \multicolumn{2}{c|}{$4$ time steps} & \multicolumn{2}{c|}{$8$ time steps} \\
 \cline{2-7}
 Level & Solution & Time & Solution  & Time & Solution & Time \\
 \hline
 $16$ & $54.923235$ & $3812.17$ & $53.634032$ & $7565.44$ & $53.512820$ & $15084.20$\\
\hline
 $17$ & $39.023559$ & $14245.68$ & $41.296346$ & $28565.06$ & $41.328455$ & $56819.84$\\
\hline
 $18$ & $39.780626$ & $55108.12$ & $38.377139$ & $109843.14$ & $38.823376$ & $220825.17$\\
\hline
 $19$ & $41.230567$ & $285159.96$ & $41.970631$ & $570320.92$ & $41.599974$ & $1140639.84 $ \\
\hline
\hline
\multicolumn{7}{|c|}{$F^{max} = 0.02$}\\
\hline
 & \multicolumn{2}{c|}{$2$ time steps} & \multicolumn{2}{c|}{$4$ time steps} & \multicolumn{2}{c|}{$8$ time steps} \\
 \cline{2-7}
 Level & Solution & Time & Solution & Time & Solution & Time \\
 \hline
 $16$ & $42.250960$ & $3812.17$ & $43.016757$ & $7565.44$ & $42.777426$ & $15084.20$\\
\hline
 $17$ & $42.116312$ & $14245.68$ & $41.406882$ & $28565.06$ & $41.842625$ & $56819.84$\\
\hline
 $18$ & $42.991274$ & $55108.12$ & $42.746717$ & $109843.14$ & $42.737987$ & $220825.17$\\
\hline
 $19$ & $42.446354$ & $285159.96$ & $42.892002$ & $570320.92$ & $42.838119$ & $1140639.84 $ \\
\hline
\end{tabular}
}
\end{center}
\caption{Sparse grid combination technique, $0.5\times 2.5$ swaption, $\sigma=0.3$, $\phi_1 =\ldots =\phi_5= 0.4$.}
\label{table:sg_d6}
\end{table}

In order to price interest rate derivatives involving more underlying forward rates using these approach, the proposed algorithm should be implemented to run on a cluster of processors (distributed memory machines). Since the communications between processors is minimal, the technique scales optimally. This extra layer of parallelism would bring also a further reduction on the previous execution times, thus allowing to stress the method with higher resolution levels.


\subsection{Numerical results with the modified sparse grid combination technique}

Our last set of numerical experiments aims at showing that the modified sparse grid technique defined by \eqref{eq:modifiedCombTechnique} is able to improve the performance (accuracy and computing time) of the standard sparse grid combination technique given by expression \eqref{eq:standardCombTechnique}, specially in moderately high dimensions.

As in the previous cases, we start with the sanity test of the pricing of the caplet with expiry $T_1$ under the classical LMM. Table \ref{table:d1:Level} gathers the behaviour of the modified technique when pricing this caplet. Firstly, we compare Table \ref{table:d1:Level} with Table \ref{table:d2_1} originated with the standard combination technique. With $\psi=1$, the modified technique is able to obtain an accuracy of $1.047853 \times 10^{-5}$ with level equal $12$ in less than a hundred seconds. In contrast, the standard sparse grid technique required a higher refinement level of $14$ and employed more than five hundred seconds to obtain a similar accuracy. Also note that with $\psi=2$ and the refinement level $10$, the modified combination technique is able to get better results, an error of $7.625043 \times 10^{-6}$ in just over $25$ seconds. Moreover, while the obtained order of convergence in space for the standard combination technique is slightly worse than two, with this modified method is almost two when $\psi=2$. The comparison with full grid method results shown in Table \ref{table:fg_sigma0} could be summarized by noting that with $\psi=2$ the modified sparse grid technique is able to obtain an error less than $1.735259 \times 10^{-6}$ in less than five hundred seconds, while the full grid approach needed almost $2.9$ hours.

\begin{table}[!htb]
\begin{center}
{\scriptsize
\begin{tabular}{|r|r|r|r|r|r|r|r|r|r|}
\hline
 & \multicolumn{4}{c|}{$\psi=1$}  & \multicolumn{4}{c|}{$\psi=2$}\\
 \cline{2-9}
 $n$ & Solution & Error & Time & \#points & Solution & Error & Time & \#points \\
 \hline
 $7$ & $6.056324$ & $2.55 \times 10^{-3}$ & $0.20$ & $2817$ & $6.057952$ & $9.25 \times 10^{-4}$ & $0.70$ & $10241$\\ %
\hline
 $8$ & $6.058978$ & $1.00 \times 10^{-4}$ & $0.44$ & $6145$ & $6.058684$ & $1.93 \times 10^{-4}$ & $2.04$ & $22529$ \\ %
\hline
 $9$ & $6.058005$ & $8.72 \times 10^{-4}$ & $1.61$ & $13313$ & $6.058822$ & $5.51 \times 10^{-5}$ & $7.08$ & $49153$\\ %
\hline
 $10$ & $6.058666$ & $2.11 \times 10^{-4}$ & $6.35$ & $28673$ & $6.058870$ & $7.62 \times 10^{-6}$ & $26.06$ & $106497$\\ %
\hline
 $11$ & $6.058822$ & $5.52 \times 10^{-5}$ & $24.62$ & $61441$ & $6.058874$ & $4.08 \times 10^{-6}$ & $105.78$ & $229377$\\ %
\hline
 $12$ & $6.058867$ & $1.04 \times 10^{-5}$ & $96.95$ & $131073$ & $6.058877$ & $9.15\times 10^{-7}$ & $478.10$ & $491521$\\ %
 \hline
 $13$ & $6.058873$ & $4.45 \times 10^{-6}$ & $473.16$ & $278529$ & $6.058877$ & $3.17 \times 10^{-7}$ & $1991.75$ & $1048577$\\ %
 \hline
 $14$ & $6.058877$ & $1.07 \times 10^{-6}$ & $1968.99$ & $589825$ & $6.058878$ & $9.91 \times 10^{-8}$ & $7868.34$ & $2228225$\\ %
 \hline 
\end{tabular}
}
\end{center}
\caption{\hspace{-1cm}Modified sparse grid combination technique, caplet with maturity $T_1$, $\sigma=0$, $256$ time steps, $F^{max}=0.04$.}
\label{table:d1:Level}
\end{table}

Finally, Tables \ref{table:d4:level} and \ref{table:d6:level} show the results for $4$ and $6$ dimensional PDEs in space, respectively. These Tables are to be compared with the corresponding Tables \ref{table:sg_d4} and \ref{table:sg_d6} generated with the standard combination technique. We observe that the higher the dimensionality of the problem the lower $\psi$ should be, otherwise the curse of dimensionality will appear soon again. Table \ref{table:d4:level} shows that with $\psi=1$ the modified method is able to obtain in just $15$ minutes a similar accuracy to the one obtained with the standard approach in more than $246$ minutes. With $\psi=2$ we observe that $5$ decimal digits are stabilized in the modified combination technique. Regarding Table \ref{table:d6:level}, with $\psi=1$ three decimal digits are stabilized already in level $13$. Also note that for the 6 dimensional PDEs in space reported in Table \ref{table:d6:level}, in our machine the modified method is not able to go further level $12$ with $\psi=2$ in a reasonable computational time due to the curse of dimensionality. Nevertheless, the accuracy recovered for level $12$ is remarkable. Finally, notice that adding points entails better performance than increasing the level of the sparse grid.

\begin{table}[!htb]
\begin{center}
{\footnotesize
\begin{tabular}{|r|rrr|rrr|}
\hline
 & \multicolumn{3}{c|}{$\psi=1$} & \multicolumn{3}{c|}{$\psi=2$} \\
 \cline{2-7}
  Level & Solution & Time & Grid points & Solution & Time & Grid points\\
\hline
 $9$ & $21.511031$ & $6.11$ & $114689$ & $21.519347$ & $107.59$ & $763905$ \\ %
\hline
 $10$ & $21.297979$ & $19.02$ & $262145$ & $21.491243$ & $328.03$ & $1765377$ \\ %
\hline
 $11$ & $21.600032$ & $60.94$ & $593921$ & $21.485699$ & $1051.35$ & $4038657$ \\ %
\hline
 $12$ & $21.595254$ & $209.27$ & $1335297$ & $21.466814$ & $3807.11$ & $9158657$ \\ %
\hline
 $13$ & $21.472738$ & $806.72$ & $2981889$ & $21.464192$ & $13921.97$ & $20611073$ \\ %
\hline
\end{tabular}
}
\end{center}
\caption{Modified sparse grids combination technique, swaption $0.5\times 1.5$, $\sigma=0.3$, $\phi_i = \phi_2 = \phi_3 =0.4$, $16$ time steps, $F^{max}=0.04$.}
\label{table:d4:level}
\end{table}

\begin{table}[!htb]
\begin{center}
{\footnotesize
\begin{tabular}{|r|rr|rr|}
\hline
 & \multicolumn{2}{c|}{$\psi=1$} & \multicolumn{2}{c|}{$\psi=2$} \\
 \cline{2-5}
  Level & Solution & Time & Solution & Time \\
\hline
 $12$ & $42.624837$ & $3255.48$ & $42.723635$ & $420696.70$ \\ %
\hline
 $13$ & $42.828046$ & $10476.64$ & $-$ & $-$ \\ %
\hline
 $14$ & $42.686665$ & $34859.98$ & $-$ & $-$ \\ %
\hline
 $15$ & $42.702808$ & $121471.65$ & $-$ & $-$ \\ %
\hline
\end{tabular}
}
\end{center}
\caption{Modified sparse grid combination technique, swaption $0.5 \times 2.5$, $\sigma=0.3$, $\phi_1 =\ldots=\phi_5= 0.4$, $4$ time steps, $F^{max} = 0.02$.}
\label{table:d6:level}
\end{table}


\section{Conclusions}
\label{sec:conclusions}

In this work we have mainly developed a new numerical \\methodology which combines high order time discretization algorithms with a sparse grids modified combination technique to solve high dimensional PDE problems arising in finance. More precisely, we have focused on the numerical solution of the PDE formulation proposed in \cite{lopezvazquez18} for pricing a large variety of interest rate derivatives, when the underlying forward rates follow a SABR-LMM model. For this purpose, we have proposed the use of high order in time AMFR-W methods, thus allowing the use of larger time steps. Moreover, a suitable splitting of the involved operators additionally contributes to the computational time reduction for a given accuracy. As the PDE problem becomes high dimensional in space when the particular interest rate derivative requires the consideration of a large number of forward rates (each one giving rise to one spatial dimension), the application of AMFR-W methods on sparse grids with combination technique turns out to be very efficient to obtain the pricing in reasonable computational times. As illustrated in the section of numerical results, parallel implementations of the algorithms based on OpenMP framework lead to a significant speed up of the computations. As indicated, an appropriate load imbalances management provides an optimal speed up, which is almost equal to the number of available computer cores. All computer implementations have been carried out from scratch. Another relevant innovative aspect comes from the suitable consideration of new homogeneous Neumann boundary conditions, instead of Dirichlet ones in \cite{lopezvazquez18}. This consideration avoids the numerical difficulties associated to the presence of boundary layers in the outflow boundaries when parameter $\beta$ is not zero, specially in the advection dominated regime. Moreover, they motivate the introduction of a modified combination technique to cope with a certain decrease in the accuracy of the standard combination technique, which mainly comes from the inaccuracy of approximations obtained with some degenerated grids included in the sparse grids combination expression. Numerical results also illustrate the advantages of the proposed modified combination technique with respect to the standard version.

Although this article focuses on the PDE formulation of the very relevant financial problem of pricing interest rate derivatives with SABR-LIBOR model, the proposed methodology can be applied to a large variety of models involving high dimension PDE formulations not only in finance but also in other disciplines in sciences and engineering. For example, in finance high dimension PDE problems related to the pricing of basket options or the computation of the XVA associated to portfolios could be considered. In computational biology, the same happens with problems related to gene networks or synthetic biology.


\appendix
\section{Appendix} \label{apend1}


\noindent {\bf Proof of Lemma \ref{lemmabij}:} Since $M_i-m_i\ge 0$, $j_i-m_i\ge 0$, for all $i$, it is clear that $\vartheta(m_1,m_2,\dots,m_N)\le \vartheta(\mathbf{j})\le \vartheta(M_1,M_2,\dots,M_N)$, for all $\mathbf{j}\in {\cal I}_N$. Trivially, we have $\vartheta(m_1,m_2,\dots,m_N)=m_1$. If we denote $Q_l:=M_l-m_l\ge 0$, with $1\le l\le N$, it is also quite simple to obtain that
\begin{align*}
&\vartheta(M_1,M_2,\dots,M_N)=M_1+\ds\sum_{l=2}^N \left( Q_l \prod_{r=1}^{l-1} (Q_r+1) \right) - (m_1-1) +(m_1-1) \\
&\qquad= \left( (Q_1+1) +\ds\sum_{l=2}^N \left( Q_l \prod_{r=1}^{l-1} (Q_r+1) \right) \right) + m_1-1 = M_T+m_1-1.
\end{align*}
On the other hand, if $\vartheta(\mathbf{j})=\vartheta(\mathbf{k})$ for two multi-indices $\mathbf{j},\mathbf{k} \in {\cal I}_N$, we have that
\begin{equation}\label{qus}
q_1=- \ds\sum_{l=2}^N \left( q_l \prod_{r=1}^{l-1} (Q_r+1) \right),\qquad q_l:=j_l-k_l,\,\forall 1\le l\le N,
\end{equation}
where the differences $q_l$ are integers that satisfy $|q_l|\le Q_l,\,1\le l\le N$. For the sake of brevity, let us suppose that $N >2$ (when $N=2$ is much simpler). From \eqref{qus}, we obtain that
\begin{equation}\label{qus2}
q_1=- \left( q_2 +\ds\sum_{l=3}^N \left( q_l \prod_{r=2}^{l-1} (Q_r+1) \right) \right) (Q_1+1).
\end{equation}
Therefore, $q_1$ is a multiple of the positive integer $(Q_1+1)$, Moreover, as $-Q_1\le q_1\le Q_1$, so necessarily $q_1$ must be zero. Therefore, because of \eqref{qus2}, we obtain
$$q_2=- \ds\sum_{l=3}^N \left( q_l \prod_{r=2}^{l-1} (Q_r+1) \right),$$
which is the same formula as in \eqref{qus}, although starting from $r=2$ instead of $r=1$. Applying a similar procedure as for $q_1$, we get that necessarily $q_2=0$. Inductively, we obtain that $q_l=0,\,i=1,2,\dots,N$, so $\mathbf{j}=\mathbf{k}$ and the map $\vartheta$ is injective. Since clearly the two sets ${\cal I}_N$ and $\{m_1,m_1+1,\dots,M_T+m_1-1\}$ have the same number $M_T$ of elements,  $\vartheta$ is a bijection. \hfill $\Box$

As a consequence of Lemma \ref{lemmabij}, for all $J\in  \{m_1,m_1+1,\dots,M_T+m_1-1\}$, there exists a unique multi-index $\mathbf{j} \in
{\cal I}_N$ given by $\mathbf{j}=\vartheta^{-1}(J)$. In practice, it is necessary to compute this inverse when we manipulate finite differences. An efficient way to calculate it is to use the \texttt{modulo}
operation, i.e., ($a \, \texttt{mod} \,n$) is  the remainder of the Euclidean division of $a$ by $n$.
\begin{lemma}\label{algor1} For every integer $J\in  \{m_1,m_1+1,\dots,M_T+m_1-1\}$, the components of the unique multi-index $\mathbf{j}=(j_1,\dots,j_N)=\vartheta^{-1}(J)\in {\cal I}_N$ can be written as $j_i=m_i+q_i$, $i=1,\dots,N$ where the integers $q_i$ satisfy
\begin{equation}\label{invj}
\begin{array}{ll}
c_1=J-m_1, & \quad q_1= c_1 \, \hbox{\rm \texttt{mod}} \, (M_1-m_1+1), \\[0.5pc]
c_i=\ds\frac{c_{i-1}-q_{i-1}}{M_{i-1}-m_{i-1}+1}, & \quad q_i= c_i \, \hbox{\rm\texttt{mod}}\, (M_i-m_i+1),\quad i=2,\dots,N. \\
\end{array}
\end{equation}
\end{lemma}

%
%
\noindent {\bf Proof of Lemma \ref{algor1}:} By using the same notation as in the proof of Lemma \ref{lemmabij}, $Q_l=M_l-m_l\ge 0,\,q_l=j_l-m_l \in \{0,1,\dots,Q_l\}$ for all $l=1,\dots,N$, equation \eqref{bij0} turns into
$$c_1:=J-m_1=q_1+ \ds\sum_{l=2}^{N} q_l \prod_{r=1}^{l-1} (Q_r+1)=q_1+ (Q_1+1) \left( q_2+ \ds\sum_{l=3}^{N} q_l \prod_{r=2}^{l-1} (Q_r+1) \right), $$
then $q_1=c_1\, \hbox{\rm \texttt{mod}} \,(Q_1+1)\in\{0,1,\dots,Q_1\}$, $j_1=q_1+m_1,$
and
$$c_2:=\ds\frac{c_1-q_1}{Q_1+1}=q_2+ \ds\sum_{l=3}^{N} q_l \prod_{r=2}^{l-1} (Q_r+1).$$
Clearly, we can apply again the \texttt{mod} operation to $c_2$, obtaining $q_2=c_2 \, \hbox{\rm \texttt{mod}} \, (Q_2+1)\in\{0,1,\dots,Q_2\}$ and $j_2=q_2+m_2$: So, if we repeat iteratively this process until the index $N$, we get  $c_N:=(c_{N-1}-q_{N-1})/(Q_{N-1}+1)=q_N$ and $j_N=q_N+m_N$. \hfill $\Box$

\begin{lemma}\label{lema:ls}
Let us consider an $s-$stage AMFR-W-method \eqref{amfrw} applied on the IVP \eqref{odeNd}-\eqref{splitbcNd}-\eqref{ffis}-Sketch 1 of dimension $M$.

For each $i=1,\dots,N$, the solution of every directional linear system of type $(I-\nu {\Delta t} {\cal A}_i)K=G$ of dimension $M$ can be obtained from the  following procedure:

For each multi-index $(j_1,\dots,j_{i-1},j_{i+1},\dots,j_N)$ of dimension $(N-1)$ with $1\le j_l\le M_l$, let us denote
\begin{equation}\label{reorder1}
J_{j_i}=\vartheta_0((j_1,\dots,j_{i-1},j_i,j_{i+1},\dots,j_N)),\qquad j_i=1,\dots,M_i,
\end{equation}
and solve the following tridiagonal system of $M_i$ equations:
\begin{equation}\label{decomp1}
\begin{array}{rl}
\left(1+2 \nu{\Delta t} \left(\delta_i\right)_{J_1}\right)  K_{J_1}  - \nu{\Delta t} \left(\delta_i\right)_{J_1} K_{J_1+E_i}&=G_{J_1}, \\[0.7pc]
- \nu{\Delta t}\left(\delta_i\right)_{J_{j_i}} K_{J_{j_i}-E_i}  + \left(1+2 \nu{\Delta t} \left(\delta_i\right)_{J_{j_i}} \right) K_{J_{j_i}}  - \nu{\Delta t} \left(\delta_i\right)_{J_{j_i}} K_{J_{j_i}+E_i}&=G_{J_{j_i}}, \\[0.5pc]   2 \,\le j_i &\le M_i-1, \\[0.7pc]
- 2\nu{\Delta t} \left(\delta_i\right)_{J_{M_i}} K_{J_{M_i}-E_i}+ \left(1+2 \nu{\Delta t}  \left(\delta_i\right)_{J_{M_i}} \right) K_{J_{M_i}}    &=G_{J_{M_i}},
\end{array}
\end{equation}
where $\left(\delta_i\right)_J=\left(d_i\right)_J /h_i^2$.

Then, the computational cost of each linear system $(I-\nu {\Delta t} {\cal A}_i)K=G$ of dimension $M$ is the cost of solving $\hat{L}_i= \prod_{k\ne i}^N M_k$ tridiagonal linear systems of dimension $M_i$.

\end{lemma}
\noindent {\bf Proof of Lemma \ref{lema:ls}:}
First, from \eqref{amfrw}, all the right hand sides $G$ of the linear systems (\ref{decomp1}) are of the form $G={\cal A}U$ for some vector $U\in \RR^L$. Therefore, for all $J=0,1,\dots,M$, we have:

If $J\in \Out$, the $J-$th row of ${\cal A}_i$ is null and $G_J=0$, so $K_J=G_J=0$.

If $J\in \Inn$,
when $J\pm E_i\in \Inn$,  ($2\le j_i\le M_i-1$), the $J-$th identity $\left((I-\nu {\Delta t} {\cal A}_i)K\right)_J=G_J$ results into
$$-\nu{\Delta t}   \left(\delta_i\right)_J K_{J-E_i} + \left(1+2 \nu{\Delta t} \left(\delta_i\right)_J \right)  K_{J}- \nu{\Delta t} \left(\delta_i\right)_J K_{J+E_i}=G_J.$$
%
It must be observed that when $(J- E_i) \in \Out$, ($j_i=1$), $K_{J- E_i}=0$, so that
$ \left(1+2 \nu{\Delta t} \left(\delta_i\right)_J \right)  K_{J}- \nu{\Delta t} \left(\delta_i\right)_J K_{J+E_i}=G_J.$
On the other hand, when $j_i=M_i$, $\vartheta_0^{-1}( J+E_i) \notin {\cal I}_N^{(0)}$, the $J$-th identity takes the form
$$-2\nu{\Delta t}   \left(\delta_i\right)_J K_{J-E_i} + \left(1+2 \nu{\Delta t} \left(\delta_i\right)_J \right)  K_{J}=G_J.$$
Then, to complete the proof it is enough to sort out these equations by taking the groups of $M_i$ equations given in \eqref{reorder1}-\eqref{decomp1}. \hfill $\Box$


In order to help the reader interested in computing the linear systems \eqref{reorder1}-\eqref{decomp1}, the procedure to solve them is presented in \cref{alg:solveLinearSystem}, when the $i-$direction and the right-hand side vector $G$ are given.

\begin{algorithm}
\caption{Procedure to solve linear systems \eqref{reorder1}-\eqref{decomp1}.}
\label{alg:solveLinearSystem}
\begin{algorithmic}
\STATE{Define a matrix $Q$ of dimension $M_i$, $K=G$}
\IF{$i=1$} \STATE {$\hat{L}_r=M_r,\quad r=2,\dots,N$ $\qquad$} \ELSIF{$i\ge 2$} \STATE{$\begin{array}{ll}
 \hat{L}_r=M_{r-1},&\quad r=2,\dots,i \\
  \hat{L}_r=M_{r},&\quad r=i+1,\dots,N
\end{array}$} \ENDIF
\STATE{$\hat{L}=\prod_{r=2}^N \hat{L}_r = \prod_{r\ne i} M_r$}
\FOR{$I=1,\dots,\hat{L}$}
 \STATE {$(k_2,\dots,k_N)=\vartheta_{N-1}^{-1}(I)$, $Q=0$}
 \STATE {$j_r=k_{r+1},\, r=1,\dots,i-1$}
 \STATE {$j_r= k_{r},\, r=i+1,\dots,N$}
 \FOR{$j_i=1,\dots,M_i$}
  \STATE {$\mathbf{j}=(j_1,\dots,j_i,\dots,j_N)$, $J=\vartheta_0(\mathbf{j})$, $R(j_i)=G(J)$}
  \STATE{$P=\left\{\begin{array}{ll} \ds\frac{\left(d_i\right)_J}{h_i^2}= \ds\frac{\alpha_i^2}{2} j_i^2 F_{N,j_N}^2, & \hbox{if } 1\le i\le N-1 \\[1pc]
\ds\frac{\left(d_N\right)_J}{h_N^2}= \ds\frac{\sigma^2}{2} j_N^2 & \hbox{if } i= N
\end{array}\right.$}
  \STATE{$Q(j_i,j_i)=-2P$}
  \IF{$j_i\ge 2$} \STATE {$Q(j_i,j_i-1)=\left\{\begin{array}{ll}       
P & \hbox{if }  j_i\le  M_i-1 \\
2P & \hbox{if }  j_i= M_i
\end{array}\right.$} \ENDIF
  \IF {$j_i\le M_i-1$} \STATE {$Q(j_i,j_i+1)=P$} \ENDIF
 \ENDFOR
 \STATE{Solve $(I_{L_i}-\nu{\Delta t} Q)X=R$}
 \FOR{$j_i=1,\dots,M_i$}
   \STATE {$\mathbf{j}=(j_1,\dots,j_i,\dots,j_N)$, $J=\vartheta_0(\mathbf{j})$, $K(J)=X(j_i)$}
 \ENDFOR
\ENDFOR
\end{algorithmic}
\end{algorithm}

As a final note, we must observe that a new bijection $\vartheta_{N-1}$ of type \eqref{bij0} is used in this algorithm to apply the re-ordering given in  \eqref{reorder1}, that is
%
$$\vartheta_{N-1}: {\cal I}_{N-1}^{(1)}=\{ \mathbf{k}=(k_2,\dots,k_N) \,\,|\,\, 1\le k_i \le \hat{L}_i,\,\,\forall i=2,\dots,N \} \longrightarrow \{1,2,\dots,\hat{L}\}, $$
with $\hat{L}=\prod_{r=2}^N \hat{L}_r$, and, for Lemma \ref{lemmabij},
$$\vartheta_{N-1} (\mathbf{k})=k_2+\ds\sum_{l=3}^N \left( (k_l-1) \prod_{r=2}^{l-1} \hat{L}_r \right),\quad \mathbf{k} \in {\cal I}_{N-1}^{(1)} .$$


\clearpage
\bibliographystyle{siamplain}
\bibliography{references}
\end{document}